\documentclass{article}
\usepackage{graphicx}
\usepackage{psfrag}
\usepackage{amssymb}
\usepackage{amsmath}
\usepackage{amsthm}
\usepackage{amscd}

\newtheorem{thm}{Theorem}[section]
\newtheorem{prop}[thm]{Proposition}
\newtheorem{cor}[thm]{Corollary}
\newtheorem{lem}[thm]{Lemma}
\newtheorem{qstn}[thm]{Question}

\newtheorem*{ec}{Ending Lamination Conjecture}
\newtheorem*{tc}{Tameness Conjecture}
\newtheorem*{thm;filling}{Theorem \ref{thm;filling}}
\newtheorem*{thm;surgery}{Theorem \ref{thm;surgery}}
\newtheorem*{twistthm}{Corollary \ref{cor;twist}}
\newtheorem*{SGL}{Theorem \ref{thm;SGL}}
\newtheorem*{thm;klarreich}{Theorem \ref{thm;klarreich}}
\newtheorem*{mainthm}{Main Theorem}

\theoremstyle{definition}

\newtheorem{exmp}[thm]{Example}

\newcommand{\A}{\mathcal{A}}
\newcommand{\T}{\mathcal{T}}
\newcommand{\TC}{\overline{\mathcal{T}}}
\newcommand{\ML}{\mathcal{ML}}
\newcommand{\MF}{\mathcal{MF}}
\newcommand{\GL}{\mathcal{GL}}
\newcommand{\PML}{\mathcal{PML}}

\newcommand{\UML}{\mathcal{UML}}
\newcommand{\EL}{\mathcal{EL}}
\newcommand{\C}{\mathcal{C}}
\newcommand{\MCG}{\mathcal{MCG}}
\newcommand{\el}{\ell}
\newcommand{\BC}{\partial_\infty{\mathcal{C}}}
\newcommand\N{{\mathbb N}}
\newcommand\Z{{\mathbb Z}}
\newcommand\R{{\mathbb R}}
\newcommand{\s}{\sigma}
\newcommand{\lam}{\lambda}
\newcommand{\mb}{\mathbf}
\newcommand{\mc}{\mathcal}

\title{The Thurston Boundary of Teichm\"uller Space and Complex of Curves}
\author{Young Deuk Kim \\ Department of Mathematics\\ Yale University\\ 
P.O. Box 208283\\New Haven, CT 06520-8283, USA 
\\(youngdeuk.kim@yale.edu)}
\date{\today}

\begin{document}
\maketitle

\begin{abstract}
Let $S$ be a closed orientable surface with genus $g\geq 2$. 
For a sequence $\s_i$ in the Teichm\"uller space of $S$, which converges to
a projective measured lamination $[\lam]$ in the Thurston boundary,
we obtain a relation between $\lam$ and the geometric limit of pants 
decompositions whose lengths are uniformly bounded by a Bers constant $L$. 
We also show that this bounded pants decomposition is related to the Gromov 
boundary of complex of curves. 

\vspace{.3cm}
\noindent
2000 Mathematics Subject Classification ; 30F60, 32G15, 57M50, 57N05. 
\end{abstract}

\section{Introduction} 

Let $S$ denote a closed orientable surface with genus $g\geq 2$, and $\T(S)$ 
the Teichm\"uller space of $S$. Although there exist some similarities between 
$\T(S)$ and a complete negatively curved space 
(see \cite{bers,FLP,krushkal,masur1}), Masur (\cite{masur}) showed that $\T(S)$
 is not negatively curved. In \cite{gromov}, Gromov introduced Gromov 
hyperbolic spaces which include complete negatively curved spaces and metric 
trees, but Masur and Wolf (\cite{MW}) showed  that $\T(S)$ is not Gromov 
hyperbolic, either (see \cite{MP1}, for another proof by McCarthy and 
Papadopoulos). 
 
\indent
In \cite{luo}, Luo classified surface theories into geometric theory, i.e. 
$\T(S)$,  algebraic theory, i.e. the mapping class group $\MCG(S)$, and 
topological theory, i.e. the complex of curves $\C(S)$. 
The complex of curves, which was introduced by Harvey in \cite{harvey}, is a 
finite dimensional simplicial complex  whose vertices are non-trivial homotopy 
classes of simple closed curves which are not boundary-parallel, and 
$k$-simplices are $k+1$ distinct vertices with disjoint representatives. 
Let $\C_0(S)$ denote the set of vertices and $\C_1(S)$ its 1-skeleton.
In \cite{MM}, Masur and Minsky defined a metric on $\C(S)$ by making each 
simplex regular Euclidean with side length 1 and taking shortest-path metric, 
and showed that $\C(S)$ is a non-proper Gromov hyperbolic space (see 
\cite{bowditch}, for a shorter proof by Bowditch). 
Every Gromov hyperbolic space has a natural boundary which is called Gromov
boundary (see \cite{BS,gromov,short}). We write $\BC(S)$ to denote the Gromov 
boundary of $\C(S)$.

\indent
Since $\T(S)$ is not Gromov hyperbolic, we can not define its Gromov boundary. 
For example, Kerckhoff(\cite{kerckhoff}) showed that the Teichm\"uller
boundary, which is the set of endpoints of geodesic rays, depends on the 
choice of base point. In \cite{thurston2}, Thurston introduced a 
compactification of $\T(S)$ with the boundary equal to the space of projective 
measured laminations $\PML(S)$, on which the action of $\MCG(S)$ extends 
continuously. Throughout this thesis, we will write $\TC(S)$ to denote this 
compactification of Thurston. See \cite{bers1,brock2} for a similar but 
different compactification by Bers. 

\indent
Thurston boundary $\PML(S)$ is the space of projective classes of measured 
laminations. We write $\ML(S)$ to denote the space of measured laminations, 
and $[\lam]$ the projective class of $\lam\in\ML(S)$. 
A measured lamination consists of a geodesic lamination and a transverse 
measure with full support on it. 
The topology on the space of geodesic laminations $\GL(S)$ is the Hausdorff 
metric topology on closed subsets. On $\ML(S)$, Thurston gave the weak-topology
induced by the measures on transverse arcs. Note that $\PML(S)$ has the 
natural quotient topology. Let $\UML(S)$ be the quotient space of $\PML(S)$ by 
forgetting measure. Although $\UML(S)$ is a subset of $\GL(S)$,
the quotient topology on $\UML(S)$ is not equal to the subspace topology.
   
\indent
In \cite{thurston2} \S 5, Thurston wrote {\em  ``... Intuitively, the 
interpretation is that a sequence of hyperbolic structures on $S$ can go to 
infinity by ``pinching'' a certain geodesic lamination $\lam$; then it 
converges to $\lam$. As a lamination is pinched toward $0$, lengths of paths 
crossing it are forced toward infinity. The ratios of these lengths determine 
the transverse invariant measure  ...'' }. 

\indent
More clearly, a sequence $\s_i\in\T(S)$ converges to $[\lam]\in\PML(S)$ in 
$\TC(S)$ if and only if for all simple closed curves $\alpha$, $\beta$ on $S$,
$${{\el_{\s_i}(\alpha)}\over {\el_{\s_i}(\beta)}}\quad\mbox{converges to}
\quad {{i(\alpha,\lam)}\over {i(\beta,\lam)}}\ ,$$ 
where $\el_{\s_i}(\alpha)$ is the length of closed $\s_i$-geodesic which is 
homotopic to $\alpha$, and $i(\alpha,\lam)$ is the intersection number of
$\alpha$ and $\lam$ which is a generalization of the geometric intersection 
number of simple closed curves (see \cite{bon,MO,papa}).  
     
\indent
The geometry of $\T(S)$ and complex of curves are well described by Minsky in 
\cite{m-geom}. They are related by the collar lemma (see \cite{buser,keen}). 
The collar lemma implies that there is a universal constant $\epsilon>0$ such 
that for any distinct $\alpha$, $\beta\in\C_0(S)$ and $\s\in\T(S)$, 
if $\el_\s(\alpha)<\epsilon$ and $\el_\s(\beta)<\epsilon$ then the two geodesic
representatives of $\alpha$ and $\beta$ are disjoint, i.e. $\alpha$ and $\beta$
are on a same simplex in the complex of curves. Therefore $\C_1(S)$ could be 
considered as the nerve of the family of regions 
$$T(\alpha)=\{\s\in\T(S)\mid\el_\s(\alpha)<\epsilon\},\quad \alpha\in\C_0(S).$$
Although $\T(S)$ is not Gromov hyperbolic, Masur and Minsky showed that $\T(S)$
is Gromov hyperbolic modulo this family of regions (see 
\cite{MM,m-geom,m-extremal}).

\indent
A lamination $\mu\in\ML(S)$ is called {\em a filling lamination} if 
$i(\mu,\mu')=0$ then
$$\mbox{support}(\mu)=\mbox{support}(\mu'),$$
for any $\mu'\in\ML(S)$.  
Minsky wrote $\EL(S)$ to denote the image of filling laminations in $\UML(S)$. 
In the celebrated proof of Thurston's ending lamination conjecture, by Brock,
Canary and Minsky, the laminations in $\EL(S)$ are appeared as ending 
laminations of Kleinian surface groups without accidental parabolics 
(see \cite{m-ending,m-ending1}). 
\begin{ec}
A hyperbolic 3-manifold with finitely generated fundamental group is uniquely 
determined by its topological type and its end invariants.
\end{ec}
\noindent
The proof of the ending lamination conjecture and the recent proof of Marden's 
tameness conjecture by Agol(\cite{agol}) give us rough picture of hyperbolic 
3-manifolds. See \cite{marden} for the original conjecture, and see \cite{CG} 
for another proof  by Calegari-Gabai. 
\begin{tc} 
A hyperbolic 3-manifold with finitely generated fundamental group is 
homeomorphic to the interior of a compact manifold with boundary.
\end{tc}
  
\indent
The Gromov boundary of $\C(S)$ is homeomorphic to the Gromov boundary of its 
1-skeleton $\C_1(S)$ because $\C(S)$ is quasi-isometric to $C_1(S)$. 
In \cite{MM}, Masur and Minsky showed that the relative hyperbolic space, which
 is roughly $\T(S)$ modulo the regions $T(\alpha)$, is quasi-isometric to 
$\C_1(S)$. 
Therefore we can expect some relation between Thurston boundary of $\T(S)$ and 
Gromov boundary of $\C(S)$ (see \cite{m-geom} for the first question of Minsky 
on this relation). In fact, Klarreich showed that $\BC(S)$ is homeomorphic to 
$\EL(S)$ (see \cite{hamenstadt} for a new proof by Hamenst\"adt). 

\begin{thm}[Klarreich \cite{klarreich}]\label{thm;klarreich}
There is a homeomorphism $$k:\BC(S)\to\EL(S)$$ such that for any sequence 
$\alpha_i$ in $\C_0(S)$, $\alpha_i$ converges to $\alpha\in\BC(S)$ if and only 
if $\alpha_i$, considered as a subset of $\UML(S)$, converges to $k(\alpha)$.
\end{thm} 

\indent
In the following theorem of Bers, $\T(\Sigma)$ stands for the Teichm\"uller 
space of $\Sigma$ with geodesic boundaries (see \cite{bers2, bers3,buser}).   
\begin{thm}[Bers \cite{bers3}]\label{thm;bers}
Let $\Sigma$ be a compact Riemann surface with genus $g\geq 2$ from which
$n$ points and $m$ disks have been removed. For any $\sigma\in\T(\Sigma)$, 
there exist $3g-3+n+m$ disjoint geodesics, which are not boundary parallel, 
whose lengths are bounded by a constant $L$ which depends only on $g,n,m$, and 
the largest length of the geodesics homotopic to the boundaries of $\Sigma$. 
\end{thm}

\indent
The constant $L$ is called {\em  a Bers constant}. Notice that for the closed 
surface $S$, for any $\s\in\T(S)$, $S$ has a pants decomposition with total 
length bounded by a constant $L$ which depends only on the genus $g$. 
In fact, Buser and Sepp\"al\"a(\cite{BS}) showed that we can choose 
$L=21g(3g-3)$. Throughout this thesis, we write $L$ to denote a fixed Bers 
constant.   

\indent
Motivated by Theorem \ref{thm;bers}, we define a function $\Phi$ on $\T(S)$ as 
follows. For $\s\in\T(S)$, let 
$$\Phi(\s)=\mbox{ a pants decomposition whose total length is bounded by L},$$ 
where all pants curves are geodesics in $\s$. Let 
$u:\PML(S)\to\UML(S)$ be the quotient map by forgetting measure. Suppose that
\begin{equation}\label{eq;basic}
\s_i\in\T(S)\mbox{ converges to }[\lam]\in\PML(S)\mbox{ in }\TC(S)
\end{equation}
and $\alpha_i$ is a pants curve in $\Phi(\s_i)$.
From Theorem \ref{thm;klarreich}, the following question is immediate.
\begin{qstn}\label{qstn1}
Suppose that $\lam$ is a filling lamination. Is the geodesic representative of 
$\alpha_i$ converging to $u([\lam])$ in $\UML(S)$? 
\end{qstn}
\noindent
This question is the motivation of all the work in this thesis. In Chapter
\ref{chap;ct}, we will solve this question positively.
Notice that if we identify $\BC(S)$ with $\EL(S)$ via the homeomorphism in 
Theorem \ref{thm;klarreich}, then we have 
\begin{thm}\label{thm;filling}
If $\lam$ is a filling lamination, then $\alpha_i$ converges to $u([\lam])$ in 
$\C(S)\cup\BC(S)$.
\end{thm}

\indent
Suppose that $\lam$ is not necessarily a filling lamination in eq. 
(\ref{eq;basic}). Is there any relation between the limit point of $\alpha_i$ 
and $u([\lam])$? 
The following simple example shows that a limit of $\alpha_i$ and $u([\lam])$
could be disjoint.
\begin{exmp}\label{exmp;disjoint}
Consider a fixed $\gamma\in\C_0(S)$. Suppose that $\s_i\in\T(S)$ converges to 
$[\gamma]\in\PML(S)$ and $\el_{\s_i}(\gamma)\to 0$, where 
we consider $\gamma$ as an element of $\ML(S)$ via the counting measure. 
Since $\el_{\s_i}(\gamma)\to 0$, there exists a pants decomposition 
$\Phi(\s_i)$ whose total length is bounded by L and
$$\gamma\in\Phi(\s_i)\quad\mbox{for all large enough }i.$$ 
Then we can find a pants curve $\alpha_i$ in $\Phi(\s_i)$ such that a limit 
point of $\alpha_i$ is disjoint from $\gamma$ in $\UML(S)$.  
\end{exmp} 

\indent
A nonempty geodesic lamination $\mu\in\GL(S)$ is called {\em  minimal} if no 
proper subset of $\mu$ is a geodesic lamination. For example, any closed 
geodesic is a minimal lamination. 
The following theorem on the structure of a geodesic lamination on $S$ works 
for any hyperbolic surface of finite type (see \cite{CEG} \S4.2, \cite{CB}
\S4 or \cite{thurston1} \S8).
\begin{thm}[Structure of Geodesic Lamination]\label{thm;SGL} 
A geodesic lamination on $S$ is the union of finitely many minimal 
sublaminations and of finitely many infinite isolated leaves whose ends spiral
along the minimal sublaminations. 
\end{thm}
\noindent
Furthermore, if $[\lam]\in\PML(S)$ then we can decompose $u([\lam])$ as a 
finite disjoint union of minimal laminations,
\begin{equation}\label{eq;decomposition}
u([\lam])=\lam_1\cup\lam_2\cup\cdots\cup\lam_m. 
\end{equation}

\begin{figure}[ht]
\begin{center}
\psfrag{F}{$F$}
\includegraphics[width=2in,height=1.5in]{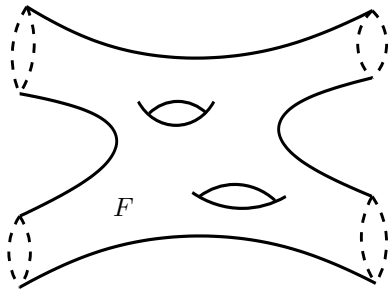}
\end{center}
\caption{An essential subsurface}
\label{fig;intro1}
\end{figure}

\indent
An essential subsurface $F$ of $S$ is a subsurface of $S$ whose boundaries are 
all homotopically non-trivial geodesics. Throughout this thesis, we assume that
all essential subsurfaces of $S$ are open, i.e. the boundaries are not 
included (see Figure \ref{fig;intro1}). 
Note that two distinct boundaries of $\overline{F}$ could be a same curve in
$S$, where $\overline{F}$ is the completion of $F$ with the path-metric in $F$.

\noindent
Suppose that $\mu\in\UML(S)$. An essential subsurface $F$ is called 
{\em filled by $\mu$}, if for any simple closed curve $\alpha$ in $F$ 
which is not parallel to a boundary of $\overline{F}$, $\alpha$ intersects 
$\mu$. 

\indent
Suppose that $\s_i\in\T(S)$ converges to $[\lam]\in\PML(S)$, and let 
$u([\lam])=\lam_1\cup\lam_2\cup\cdots\cup\lam_m$ be the decomposition into 
minimal sublaminations. The main theorem of this thesis is  

\begin{thm}[Main Theorem]\label{mainthm}
If $\Phi(\s_i)$ converges to a geodesic lamination $\nu$ in Hausdorff metric 
topology, then $u([\lam])\subset\nu$.
\end{thm}

\indent
There is no minimal lamination which fills a pair of pants (see \cite{penner}
\S 2.6), therefore in eq. (\ref{eq;decomposition}), if $\lam_j$ is not a simple
closed curve then the essential subsurface filled by $\lam_j$ is at least a 
1-holed torus or 4-holed sphere. 
For a closed annulus $Y$, Minsky defined the arc complex $\A(Y)$. 
In this complex vertices are essential homotopy classes, rel endpoints, of 
properly embedded arcs, and simplices are sets of vertices with representatives
with disjoint interiors. Note that the endpoints are {\em not} allowed to move 
in the boundary. Let $\A_0(Y)$ be the set of vertices and $\A_1(Y)$ the 
1-skeleton, and give shortest-path metrics to $\A(Y)$ and $\A_1(Y)$ as in 
$\C(S)$. For $a,b\in\A_0(Y)$, we write $a\cdot b$ to denote the {\em algebraic}
intersection number. 

\indent
Suppose that in eq. (\ref{eq;decomposition}), $\lam_1$ is a simple closed 
curve. Consider an annular covering $Y$ of $S$ in which a neighborhood of 
$\lam_1$ lifts homeomorphically. The following corollary follows from the
main theorem. 

\begin{cor}\label{cor;twist}
Suppose that $\alpha_i$ meets $\lam_1$ for all $i$. Then $|a_1\cdot a_i|$ 
approaches to $\infty$, where $a_1\in\A_0(Y)$ is a lift of $\alpha_1$ in $Y$ 
and so is $a_i$.
\end{cor}

\indent
Suppose that in eq. (\ref{eq;decomposition}), $\lam_2$ is not a simple closed 
curve. Suppose that $F$ is an essential subsurface of $S$ which is filled by 
$\lam_2$. Notice that there exists a pants curve $\alpha_i$ in $\Phi(\s_i)$ 
such that $\alpha_i\cap F\neq\emptyset$ for all $i$. Let $\beta_i$ be a 
component of $\alpha_i\cap F$. We can define a simple closed curve 
$\Tilde{\beta}_i$ from $\beta_i$ canonically (see Chapter \ref{chap;ct}). 
We will prove the following theorem which is a generalization of Theorem 
\ref{thm;filling}.

\begin{thm}\label{thm;surgery}
The geodesic representative of $\Tilde{\beta}_i$ converges to $\lam_2$ in 
$\UML(F)$.
\end{thm}

\indent
In Chapter \ref{chap;tbdry} and Chapter \ref{chap;gbdry}, we study some 
preliminaries. We prove the main theorem in Chapter \ref{chap;pfmainthm},
and we prove Corollary \ref{cor;twist} and Theorem \ref{thm;surgery} in 
Chapter \ref{chap;ct}. 

\vspace{.5cm}
\textbf{\large{Acknowledgements.}}
This paper is a part of Ph.D thesis of the author (May 2005, Stony Brook 
university). 
The author would like to express his deepest gratitude to his advisor 
Yair Minsky for suggesting this thesis topic, and for his patient guidance 
and continuous encouragement.

\section{The Thurston Boundary of Teichm\"uller Space}\label{chap;tbdry} 

Let $S$ be a closed oriented surface of genus $g\geq 2$. In this chapter we 
study the Teichm\"uller space $\T(S)$ and its Thurston boundary $\PML(S)$.
See \cite{IT} and \cite{thurston5} as references. 

\subsection{Teichm\"uller space and its boundaries}\label{sec;teich}

A conformal structure $\s$ on $S$ is determined by an atlas of coordinate 
neighborhoods $(U_\alpha,z_\alpha)$, where $\{U_\alpha\}$ is an open cover of 
$M$, and $z_\alpha:U_\alpha\to\mathbb{C}$ has the property that 
$z_\alpha\circ z_\beta^{-1}$ is analytic whenever defined.
Teichm\"uller space $\T(S)$ is the space of conformal structures
on $S$, where two structures are considered to be equivalent if there
is a conformal map between them isotopic to the identity. 
Recall that two diffeomorphisms $f$ and $g$ on $S$ are called {\em isotopic}
if there exists a diffeomorphism $H(x,t)=(h_t(x),t):S\times[0,1]\to S
\times[0,1]$ such that $h_0(x)=f(x)$ and $h_1(x)=g(x)$ for all $x\in S$. 

\indent
By the uniformization theorem, $\T(S)$ can be considered as the space of 
complete hyperbolic metrics on $S$ with finite area. The area of 
hyperbolic surface $(S,\s)$ is equal to $-2\pi\chi(S)$, where $\chi(S)=2(1-g)$ 
is the Euler characteristic of $S$ (see \cite{BP} for details).

\indent
As usual, two closed curves $\alpha,\beta:[0,1]\to S$ are called 
{\em free homotopic} if there exists a continuous mapping 
$F:[0,1]\times[0,1]\to S$ such that 
$$F(t,0)=\alpha(t),\ F(t,1)=\beta(t)\quad\mbox{and}\quad F(0,s)=F(1,s)$$
for all $t,s\in [0,1]$. The following well-known lemma will be used frequently 
in this thesis (see \cite{BP} \S B.4 for a proof).
\begin{lem}\label{lem;gloop}
Suppose that $(S,\s)$ is a hyperbolic surface. Then each free homotopy class of
closed curves contains a unique geodesic representative.  
\end{lem}
\noindent
Two closed curves $\alpha,\beta:[0,1]\to S$ are called {\em isotopic} if there 
exists a diffeomorphism $G:S\times[0,1]\to S\times[0,1]$ such that for all 
$x\in S$ and $t,s\in [0,1]$,
$$G(x,t)=(g_t(x),t),\quad
g_0(x)=x\quad\mbox{and}\quad g_1(\alpha(s))=\beta(s).$$
Free homotopic simple curves on a connected surface are isotopic, too.

\indent
There is well-known classification of $\mathit{Isom}{\,\mb{H}}^n$ into 
elliptic, parabolic and hyperbolic isometries. 
We write $\overline{\mb{H}}^2$ to denote the compactification of $\mb{H}^2$ 
by the circle at infinity. The following two lemmas will be useful in Section
\ref{sec;tnumber} (see \cite{BP} \S B.4 for a proof).
\begin{lem}\label{lem;hisometry}
Every non-trivial elements of $\pi_1(S)$ are hyperbolic isometry.
\end{lem}   
\begin{lem}\label{lem;lift}
Suppose that $\alpha$ is a non-trivial simple closed geodesic in $S$. Then
any two lifts of $\alpha$ into $\mb{H}^2$ can not meet in the whole
$\overline{\mb{H}}^2$.  
\end{lem}   
\noindent
The following lemma will be useful in Chapter \ref{chap;pfmainthm} 
(see \cite{BP} \S B.4 for a proof).
\begin{lem}\label{lem;disjoint}
Suppose that $\alpha$ and $\beta$ are non-intersecting, non-isotopic and 
non-trivial simple closed curves in a hyperbolic surface $(S,\s)$. Then the 
geodesic representatives of $\alpha$ and $\beta$ are non-intersecting.  
In particular, if $\alpha$ is in a subsurface $F$ of $S$, then the geodesic
representative of $\alpha$ is in $F$, too.
\end{lem}   

\indent
The standard boundary $S^{n-1}_\infty$ of $\mb{H}^n$ is defined by the 
equivalence classes of geodesic rays (see \cite{BP} \S A.5). 
The Teichm\"uller boundary of $\T(S)$ is defined in the same way. 
But there are distinct geodesic rays from a point in $\T(S)$ that always remain
within a bounded distance of each other (see \cite{masur1}). 

\indent
The mapping class group $\MCG(S)$ is the group of orientation preserving 
diffeomorphisms of $S$ modulo those which are isotopic to identity, i.e.
$$\MCG(S)=\mbox{Diff}^+(S)/\mbox{Diff}^0(S).$$
The mapping class group acts by isometries on $\T(S)$ and its quotient space is
called {\em moduli space} of $S$. See \cite{maskit3} for a picture of 
moduli space. Kerckhoff proved in his thesis that the action of $\MCG(S)$ does 
not in general extend to Teichm\"uller boundary.
 
\indent
Let $\mathcal{V}(S)$ be the set of representations of $\pi_1(S)$ into 
$PSL(2,\mathbb{C})$ up to conjugacy with compact-open topology. The product
$\T(S)\times\T(S)$ can be identified with an open subset of $\mathcal{V}(S)$, 
consisting of faithful representations whose images are quasi-Fuchsian groups, 
by Bers simultaneous uniformization. If we fix the first factor, then we get
a holomorphic embedding of $\T(S)$ into $\mathcal{V}(S)$ which is called
{\em a Bers slices}. Although this embedding depends on the fixed first factor,
there is a biholomorphic mapping between any two slices. 
The closure of this slice is compact, and is called {\em a Bers 
compactification} (see \cite{bers1,mcmullen}). 
But Kerckhoff and Thurston showed that (see Theorem 1 and Theorem 2 of 
\cite{KT})
\begin{enumerate}
\item
For each genus $g\geq 2$, there are Bers slices for which the canonical 
homeomorphisms do not extend to homeomorphisms on their compactifications. 
\item
For $g=2$, there is a Bers slice for which the action of the mapping class 
group does not extend continuously to its compactification.
\end{enumerate}  

\indent
In the famous 1976 preprint, which is published in \cite{thurston5} later, 
Thurston introduced the space of projective measured laminations on $S$, which 
will be denoted by $\PML(S)$, and a compactification of $\T(S)$ whose boundary 
is equal to $\PML(S)$. 
Thurston boundary $\PML(S)$ is a natural boundary of $\T(S)$, in the sense
that the action of mapping class group extends continuously to the Thurston 
compactification $\TC(S)=\T(S)\cup\PML(S)$. Masur(\cite{masur1}) showed that 
Teichm\"uller boundary and Thurston boundary are same almost everywhere, but 
not everywhere.

\indent
With this compactifcation, Thurston classified surface diffeomorphisms as 
periodic, reducible or pseudo-Anosov, which is a generalization of the 
well-known classification of elements of $SL(2,\Z)$. 
Thurston boundary was also used by Kerckhoff to solve the Nielsen realization 
problem, i.e. {\em every finite subgroup of $\MCG(S)$ can be realized as a 
group of isometries of some hyperbolic structure on $S$} 
(see \cite{kerckhoff2} for the proof).  

\subsection{Measured laminations}

In this section we study measured laminations. See \cite{bon4,CEG,CB,hatcher,
levitt} and \cite{thurston1} as references. 
Consider a fixed hyperbolic structure $\s$ on $S$. A {\em geodesic lamination} 
$\mu$ is a closed subset of $S$, which is a disjoint union of simple geodesics 
which are called leaves of $\mu$. The leaves of a geodesic 
lamination are complete, i.e. each leaf is either closed or has infinite 
length in both of its ends, and a geodesic lamination is determined by its 
support, i.e. a geodesic lamination is a union of geodesics in just one way.  
Using $S^1_\infty$, a geodesic lamination on $(S,\s)$ can be naturally related 
to a geodesic lamination on $(S,\s')$ for any $\s'\in\T(S)$. 
We write $\GL(S)$ to denote the space of geodesic laminations on $S$, which is
equipped with the Hausdorff metric on closed subsets.
Note that $\GL(S)$ is compact and therefore, in particular, every infinite 
sequence of nontrivial simple closed geodesics has a convergent subsequence.

\indent
For an arbitrary topological space $X$, the Chabauty topology on the set of 
closed subsets of $X$ has the following sub-bases.
\begin{enumerate}
\item[(i)]
$O_1(K)=\{A\mid A\cap K=\emptyset\}$ where $K$ is compact.
\item[(ii)]
$O_2(U)=\{A\mid A\cap U\neq\emptyset\}$ where $U$ is open.
\end{enumerate}

\noindent
If $X$ is compact and metrizable, in particular for $S$, the Chabauty topology 
agrees with the topology induced by the Hausdorff metric. The following lemma 
will turn out to be useful (see \cite{CEG} \S 3.1). 
\begin{lem}[Geometric Convergence]\label{lem;gc}
Suppose that $X$ is a locally compact metric space. A sequence $A_n$ of closed 
subsets of $X$ converges to a closed subset $A$ in Chabauty topology if and
only if 
\begin{enumerate}
\item[(i)]
If $x_{n_k}\in A_{n_k}$ converges to $x\in X$ then $x\in A$.
\item[(ii)]
If $x\in A$, then there exists a sequence $x_n\in A_n$ which converges to $x$. 
\end{enumerate}
\end{lem}

\indent
In $\mb{H}^2$, a geodesic is determined by an element of the open M\"obius band
$$M=(S^1_\infty\times S^1_\infty-\Delta)/\mb{Z}_2,$$
where $\Delta=\{(x,x)\}$ is the diagonal and $\mb{Z}_2$ acts by interchanging 
coordinates. A geodesic in $\mb{H}^2$ projects to a simple geodesic on $S$ if 
and only if the covering translates of its pairs of end points never strictly 
separate each other. 
Notice that a geodesic lamination could be considered as a closed subset of 
$M$. The Chabauty topology on $\GL(S)$ as closed subsets of $M$ is equivalent 
to the Chabauty topology on $\GL(S)$ as closed subsets of $\mb{H}^2$.
Therefore
\begin{lem}\label{lem;HCtopology}
If $\mu_i$ converges to $\mu$ in $\GL(S)$ with Hausdorff metric topology, then
for any geodesic $\el\subset\mu$, there exist geodesics $\el_i\subset\mu_i$ 
which converge to $\el$.
\end{lem}

\indent
A geodesic lamination is called {\em maximal} if each complementary 
region is isometric to an ideal triangle. A nonempty geodesic lamination is 
called {\em minimal} if no proper subset is a geodesic lamination. 
For example, any simple closed geodesic is a minimal lamination.
The following lemma is about the structure of minimal laminations 
(see \cite{CEG} \S 4.2 for a proof).   
\begin{lem}[Structure of minimal lamination]
If $\mu$ is a minimal lamination then either $\mu$ is a single geodesic or
consists of uncountable leaves. 
\end{lem}

\noindent
The following theorem is about the structure of geodesic laminations.   
\begin{SGL}[Structure of Geodesic Lamination] 
A geodesic lamination on $S$ is the union of finitely many minimal 
sublaminations and of finitely many infinite isolated leaves whose ends spiral
along the minimal sublaminations. 
\end{SGL}

\begin{figure}[ht]
\begin{center}
\psfrag{a}{$\alpha$}
\psfrag{b}{$\beta$}
\includegraphics[width=3in,height=1.5in]{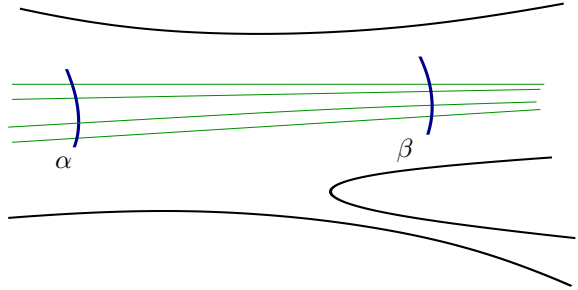}
\end{center}
\caption{$\alpha$ and $\beta$ have same measure}
\label{fig;tb1}
\end{figure}

\indent
A transverse measure on a geodesic lamination $\mu$ is a rule, which assigns to
each transverse arc $\alpha$ a measure that is supported on $\mu\cap\alpha$, 
which is invariant under a map from $\alpha$ to another arc $\beta$ if it 
takes each point of intersection of $\alpha$ with a leaf of $\mu$ to a point 
of intersection $\beta$ with the same leaf (see Figure \ref{fig;tb1}).
A measured lamination on $S$ is a geodesic lamination $\mu$ with a transverse 
measure of full support, i.e. if $\alpha\cap\mu\neq\emptyset$ then $\alpha$ 
has nonzero measure for any transverse arc $\alpha$. For example, a simple 
closed geodesic equipped with counting measure is a measured lamination. 
We write $\ML(S)$ to denote the space of measured laminations on $S$.
There is a natural action of $\R^+$ on $\ML(S)$. Suppose that $r>0$. 
The measured lamination $r\mu$ has the same geodesic lamination as $\mu$ with 
the transverse measure scaled by $r$. 
We write $\PML(S)$ to denote the set of equivalence classes of projective 
measured laminations.

\indent
The support of a measured lamination has no infinite isolated leaves,  
therefore by Theorem \ref{thm;SGL}, it is a finite disjoint union of minimal 
sublaminations. 

\indent
For a transverse arc $\alpha$ to a measured lamination $\mu$, let $\theta$ be 
the angle between leaves of $\mu$ and $\alpha$, measured counterclockwise from
$\alpha$ to $\mu$. The {\em total angle} of $\alpha$ is defined by
$$\theta(\alpha,\mu)=\int_\alpha\theta\, d\mu.$$
Thurston gave $\ML(S)$ a topology with following basis, 
\begin{eqnarray*}
&&\mathcal{B}(\mu,\alpha_1,\cdots,\alpha_n,\epsilon)=\{\nu\in\ML(S):\\
&&\qquad\qquad\qquad |\left(\mu(\alpha_k),\theta(\alpha_k,\mu)\right)-
\left(\nu(\alpha_k),\theta(\alpha_k,\nu)\right)|<\epsilon,\ k=1,\cdots, n\},
\end{eqnarray*} 
where $\{\alpha_k\}$ is a finite set of transverse arcs to $\mu$, and 
$\epsilon>0$. The following theorem was proved by Thurston.  
\begin{thm}[Thurston]
\begin{itemize}
\item[(1)]
$\ML(S)$ is homeomorphic to the open ball $B^{6g-6}$ and $\PML(S)$ is 
homeomorphic to the sphere $S^{6g-7}$.
\item[(2)]
$\R\times\C_0(S)$ is dense in $\ML(S)$, and $C_0(S)$ is dense in $\PML(S)$. 
\end{itemize}
\end{thm}   

\indent
Suppose that $\alpha,\beta\in\C_0(S)$. The {\em geometric intersection number} 
$i(\alpha,\beta)$ is the minimal number of intersections of any two their 
representatives.

\indent
For a transverse arc $\alpha$ to $\mu\in\ML(S)$, we write $\int_\alpha\, d\mu$ 
to denote integration of the transverse measure over $\alpha$. 
For a simple closed curve $\gamma$, let 
$$i(\mu,\gamma)=\inf_{\gamma'}\int_{\gamma'}d\mu,$$ 
where the infimum is taken over all the simple closed curves $\gamma'$ which is
homotopic to $\gamma$. For a general transverse arc $\alpha$, we define 
$$i(\mu,\alpha)=\inf_{\alpha'}\int_{\alpha'}d\mu,$$   
where the infimum is taken over all the arcs $\alpha'$ which is homotopic to 
$\alpha$ with endpoints fixed. Note that, in both cases, the infimum is 
realized by the unique geodesic in the corresponding homotopy class. 

\indent
Suppose that $\mu\in\ML(S)$, $\gamma\in\C_0(S)$ and $r\in\R^+$, let
$i(\mu,r\gamma)=ri(\mu,\gamma)$. 
The following theorem of Thurston is useful in Chapter \ref{chap;ct} 
(see \cite{thurston4} for a proof).   

\begin{thm}[Continuity of the intersection number]
The intersection number $i$ extends to a continuous symmetric function on 
$\ML(S)\times\ML(S)$.   
\end{thm}   

\subsection{Topology of $\TC(S)$}

In this section we study the topology of the Thurston compactification 
$\TC(S)$. References are \cite {FLP,papa} and \cite{thurston4}. The topology 
on $\TC(S)=\T(S)\cup\PML(S)$ is determined by the following two properties.
\begin{enumerate}
\item[\textbf{(P1)}]
$\T(S)$ is open in $\TC(S)$.
\item[\textbf{(P2)}]
$\s_i\in\T(S)$ converge to $[\lam]\in\PML(S)$ if and only if, for all simple 
closed curves $\alpha$, $\beta$ on $S$ with $i(\beta,\lam)\neq 0$, 
$${{\el_{\s_i}(\alpha)}\over {\el_{\s_i}(\beta)}}\quad\mbox{converges to}
\quad {{i(\alpha,\lam)}\over {i(\beta,\lam)}}\ ,$$ 
where $\el_{\s_i}(\alpha)$ is the length of closed $\s_i$-geodesic which is 
homotopic to $\alpha$.     
\end{enumerate}

\indent
The following trivial lemma will turn out to be useful (see \cite{MP} for a
proof). 
\begin{lem}
For any infinite sequence of distinct simple closed curves in $S$, there is a
subsequence $\alpha_i$ and $c_i>0$ such that 
$$c_i\to 0\quad\mbox{and}\quad c_i\alpha_i\to\mu\quad\mbox{for some }
\mu\in\ML(S)-\{0\}.$$     
\end{lem}

\indent
Suppose that $\mu_i\in\ML(S)$. We write $\mu_i\to\infty$ to denote that there 
exists $\alpha\in\C_0(S)$ such that $i(\alpha,\mu_i)$ converges to $\infty$.
The following theorem is the most useful theorem in this thesis.

\begin{thm}[Theorem 2.2 of \cite{thurston4}]\label{thm;useful}
A sequence $\s_i\in\T(S)$ converges to $[\lam]\in\PML(S)$ if and only if 
there is a sequence $\mu_i\in\ML(S)$ converging projectively to $\lam$ such 
that $\mu_i\to\infty$ and $\el_{\s_1}(\mu_i)\to\infty$ but $\el_{\s_i}(\mu_i)$ 
remains bounded, and for all $\nu\in\ML(S)$, there exists a constant $C>0$ 
such that 
$$i(\nu,\mu_i)\leq \el_{\s_i}(\nu)\leq i(\nu,\mu_i)+C \el_{\s_1}(\nu).$$
\end{thm}

\noindent
In particular, for each $\gamma\in\C_0(S)$, there exists a constant $\Gamma>0$,
which does not depend on $i$ such that 
\begin{equation}\label{eq;papa}
i(\gamma,\mu_i)\leq \el_{\s_i}(\gamma)\leq i(\gamma,\mu_i)+\Gamma.
\end{equation}

\indent
The measured laminations $\mu_i$ were constructed in \cite{FLP}, \cite{papa} 
and \cite{thurston3} \S 9. In these constructions, the associated measured 
foliations are constructed first, and the measured laminations $\mu_i$ are 
induced later. 
In the following paragraphs, we study measured foliations and the construction 
of $\mu_i$ following Papadopoulos (see \cite{papa} for details). 
Another good French reference for this construction is $\cite{FLP}$. 

\indent
A measured foliation $F$ on $S$ is a foliation with finite number of 
singularities equipped with a invariant transverse measure, i.e. $F$ is 
determined by a finite number of points $p_k\in S$ and an atlas of coordinate 
neighborhoods $$(x_i,y_i):U_i\to\R^2$$ on the complement of $\{p_k\}$ such that
$x_j=f_{ij}(x_i,y_i)$ and $y_j=\pm y_i+C$ for any overlapping coordinate 
neighborhoods $(x_j,y_j)$, where $C$ is a constant and the transverse measure
is $dy$. The singularities have $p$-pronged saddles with $p\geq 3$.   
Suppose that $F$ is a measured foliation on $S$ and $\alpha\in\C_0(S)$. Let
$$i(F,\alpha)=\inf_{\alpha'}\int_{\alpha'}|dy|,$$
where the infimum is taken over all the representatives $\alpha'$ in the class 
$\alpha$. Two measured foliations $F$ and $G$ are called {\em equivalent} if 
$i(F,\alpha)=i(G,\alpha)$ for all $\alpha\in\C_0(S)$. We write $\MF(S)$ to 
denote the set of equivalence classes of measured foliations.
Measured foliations and measured laminations are related by the following 
theorem (see \cite{levitt} for the proof).

\begin{thm}\label{thm;FL}
There is a homeomorphism $h:\MF(S)\to\ML(S)$ which is identity on 
$\R\times\C_0(S)$ and preserves the intersection number.
\end{thm} 

\begin{figure}[ht]
\begin{center}
\includegraphics[width=2in,height=2in]{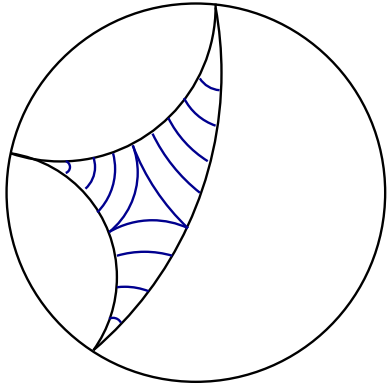}
\end{center}
\caption{The partial foliation}
\label{fig;tb2}
\end{figure}

\indent
Suppose that $\mu$ is a maximal geodesic lamination on $S$ and $\s\in\T(S)$. 
Thurston constructed a measured foliation $F_\mu(\s)$, which is called 
{\em horocyclic foliation} as follows.
Since $\mu$ is maximal, the complementary components of $\mu$ are all isometric
to ideal triangles. In each of these components, define a partial foliation,
i.e. a foliation whose support is subsurface, whose leaves are intersections of
the triangle and horocycle centered at vertices of the triangle
(see Figure \ref{fig;tb2}).

\noindent
Notice that the horocycle meets the triangle with right angles, and the 
non-foliated region is equal to a little triangle whose edges are subarcs of
horocycles which meet tangentially at their endpoints. These partial foliations
in the ideal triangles fit together on the surface and define a partial 
foliation on $S$. The transverse measure on this partial foliation is uniquely 
determined by the fact that on the leaves of $\mu$ this transverse measure is 
equal to the hyperbolic distance.  

\indent
Suppose that $\mu$ is a maximal geodesic lamination on $S$. Let $\MF(\mu)$ be
the subset of $\MF(S)$ consisting of equivalence classes which has 
representative transverse to $\mu$. In \cite{thurston3} \S 9, Thurston showed 
\begin{thm}
The map $\phi_\mu:\T(S)\to\MF(\mu)$ such that $\phi_\mu(\s)=F_\mu(\s)$
is a homeomorphism.
\end{thm}

\noindent  
Suppose that $\s_i\in\T(S)$ converges to $[\lam]\in\PML(S)$ in $\TC(S)$. 
Then the measured lamination $\mu_i$, obtained from $F_\mu(\s_i)$ by 
Theorem \ref{thm;FL}, is the lamination in Theorem \ref{thm;useful}.

\section{Complex of Curves}\label{chap;gbdry} 

In this chapter we study the complex of curves and its Gromov boundary.
All theories in this chapter work for any orientable surface of finite type, 
i.e. surface with genus $g$ and $n$ punctures. Throughout this chapter, 
we write $\Sigma=\Sigma_{g,n}$ to denote an orientable surface with genus $g$ 
and $n$ punctures. As before, we write $\MCG(\Sigma)$ to denote the mapping 
class group of $\Sigma$, and $\T(\Sigma)$ the Teichm\"uller space of $\Sigma$. 
 
\subsection{Complex of curves}\label{sec;complex of curves}

In \cite{harvey}, Harvey introduced the complex of curve $\C(\Sigma)$ to study 
the action of $\MCG(\Sigma)$ at the infinity of $\T(\Sigma)$. This complex 
encodes the asymptotic geometry of Teichm\"uller space, similarly as the 
Tits buildings for symmetric spaces. Let $\mc{S}(\Sigma)$ be the set of isotopy
classes of essential, unoriented, non-boundary parallel simple closed curves in
$\Sigma$. 
The vertices of $\C(\Sigma)$ are elements of $\mc{S}(\Sigma)$, i.e.
$\C_0(\Sigma)=\mc{S}(\Sigma)$, and the $k$-simplices of $\C(\Sigma)$ are  
subsets $\{\alpha_1,\cdots,\alpha_{k+1}\}$ of $\mc{S}(\Sigma)$ with mutually 
disjoint representatives. Notice that $\C(\Sigma)$ is empty if $g=0$ and 
$n\leq 3$. The maximal dimension of simplices is called the {\em dimension} of 
$\C(\Sigma)$, and it is equal to $3g+n-4$.
Masur and Minsky defined a metric on $\C(\Sigma)$ by making each simplex 
regular Euclidean with side length 1 and taking shortest-path metric. 
Let $d_\C$ denote this metric on $\C(\Sigma)$. 

\indent
The mapping class group acts on $\C(\Sigma)$ and Ivanov proved that, if 
$g\geq 2$ then all automorphisms of  $\C(\Sigma)$ are given by elements of
$\MCG(\Sigma)$ (see \cite{ivanov} for the proof, and \cite{korkmaz} 
for the related work of Korkmaz). Luo(\cite{luo3}) generalized this result 
by showing that, if $3g+n-4\geq 1$ and $(g,n)\neq (1,2)$, then all 
automorphisms of $\C(\Sigma)$  are given by elements of $\MCG(\Sigma)$. 
In \cite{harer}, Harer showed that $\C(\Sigma_{g,n})$ is homotopic to a 
wedge of spheres of dimension $r$, where 
\begin{eqnarray*}
r=\left\{
\begin{array}{lll}
2g+n-3 &\mbox{ if }& g>0\mbox{ and }n>0\\ 
2g-2 &\mbox{ if }& n=0\\ 
n-4 &\mbox{ if }& g=0.  
\end{array}
\right.
\end{eqnarray*}
            
\indent
If $\Sigma$ is a torus, once-punctured torus or 4-times punctured sphere, then 
any two essential simple closed curves intersects, i.e. there is no edge in 
$\C(\Sigma)$.
Masur and Minsky introduced a new definition for these cases so that it has 
edges. In this definition, $\{\alpha,\beta\}$ is an edge if $\alpha\neq\beta$, 
and $\alpha$ and $\beta$ have the lowest possible intersection number. 
For the tori this is 1, and for 4-holed sphere this is 2 (see \cite{HT} for 
Hatcher-Thurston complex which is related to this definition).

\subsection{Relative twist number in annulus complex}\label{sec;tnumber}

In this section we study the relative twist number in annulus complex, which 
was introduced by Minsky in \cite{FLM}, of two simple closed curves around a 
fixed simple closed curve. See \cite{FLM} \S2.1, \cite{MM2} \S2.4 and 
\cite{m-ending} \S4 as references. 
An annular domain in $\Sigma$ is an annulus with incompressible boundary.
A complex, which is called an {\em annulus complex}, is defined for such annuli
to keep track of Dehn twisting around their cores.

\indent
Consider an oriented annulus $Y=S^1\times [0,1]$. We write $\A_0(Y)$ to denote 
the set of arcs joining $S^1\times\{0\}$ to $S^1\times\{1\}$, up to homotopy 
with endpoints fixed. As in $\C(\Sigma)$, we put an edge between any 
two elements of $\A_0(Y)$ which have representatives with disjoint interiors, 
and define the annular complex $\A(Y)$ as for the complex of curves. 
We also make $\A(Y)$ a metric space with edge length 1 as in the curve complex.
Let $d_Y$ denote the path-metric.

\vspace{.5cm}  
\begin{figure}[ht]
\begin{center}
\psfrag{a}{$a$}
\psfrag{b}{$b$}
\includegraphics[width=2.5in,height=1in]{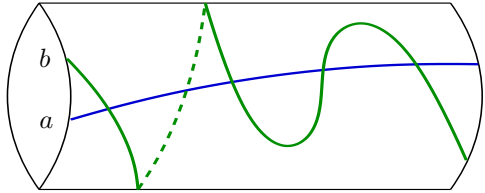}
\end{center}
\caption{$|a\cdot b|=2$}
\label{fig;cc1}
\end{figure}

\indent
Suppose that $a,b\in\A_0(Y)$ and they do not share any endpoints. Notice that 
$a$ and $b$ inherit orientations from the orientation of $[0,1]$. Therefore we 
can define the {\em algebraic intersection number} $a\cdot b$ (for example, 
see Figure \ref{fig;cc1}). Let $a\cdot a=0$.
Consider a lift of $a\in\A_0(Y)$ to the covering space 
$\widetilde{Y}=\R\times[0,1]$ which has endpoints $(a_0,0)$ and $(a_1,1)$. 
Notice that these endpoints are determined by $a$, up to $\Z$, and 
$$a\cdot b=\lfloor b_1-a_1\rfloor-\lfloor b_0-a_0\rfloor, $$
where $\lfloor x\rfloor$ denotes the largest integer less than or equal to $x$.
It follows that 
\begin{equation}\label{eq;dot}
a\cdot c=a\cdot b+b\cdot c+\Delta\quad\mbox{with }\Delta\in\{0,1,-1\}
\end{equation}
for all $a,\ b,\ c\in\A_0(Y)$ such that the intersection numbers are defined. 
With an inductive argument, we can also check  
\begin{equation}\label{eq;dotdist}
d_Y(a,b)=1+|a\cdot b|\quad\mbox{for all distinct }a,b\in\A_0(Y). 
\end{equation}
Fix $a\in\A_0(Y)$. From eq. (\ref{eq;dot}) and eq. (\ref{eq;dotdist}), we can 
show that the map $f:\A_0(Y)\to\Z$ with $f(b)=a\cdot b$ is a 
quasi-isometry. Thus $\A(Y)$ is quasi-isometric to $\Z$. 

\indent
For a fixed finite generating set of $\MCG(\Sigma)$, let
$||\cdot ||$ be the minimal word length with respect to these generators. 
Masur and Minsky introduced the relative twist number, and Farb, Lubotzky and 
Minsky proved that every Dehn twist has linear growth in $\MCG(\Sigma)$.
    
\begin{thm}[Theorem 1.1 of \cite{FLM}]\label{thm;FLM}
For all Dehn twist $t$, there exists a constant $c>0$ such that 
$||t^m||\geq c|m|$ for all $m$.
\end{thm}

\begin{figure}[ht]
\begin{center}
\psfrag{a}{$a$}
\psfrag{b}{$b$}
\psfrag{Y}{$Y_\alpha$}
\psfrag{H}{$\overline{\mb{H}}^2$}
\includegraphics[width=2in,height=2.5in]{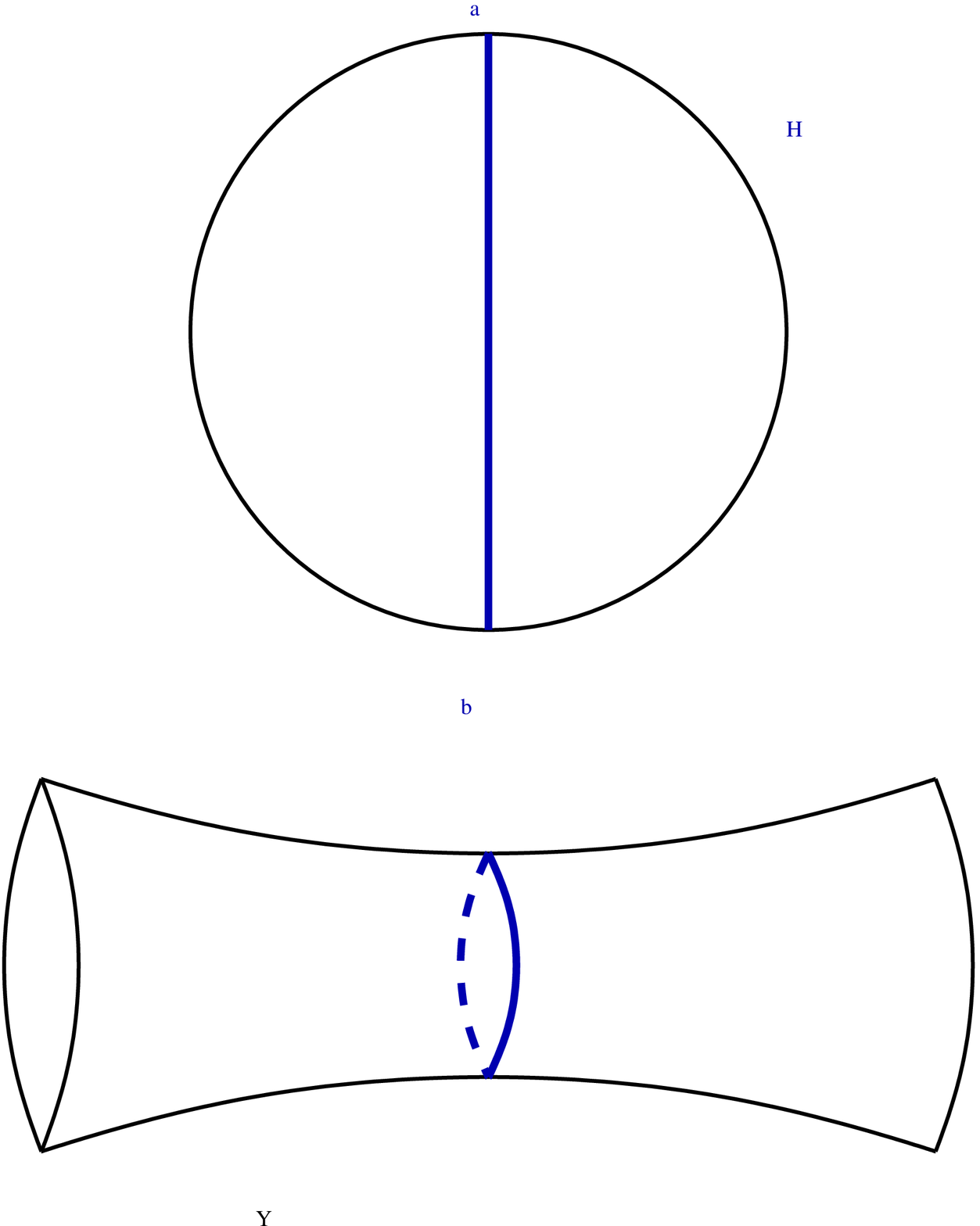}
\end{center}
\caption{$Y_\alpha$}
\label{fig;cc2}
\end{figure}

\indent 
The relative twist number is defined as follows. For a fixed essential simple 
closed curve $\alpha$ in $\Sigma$, let $g_\alpha$ be an isometry of $\mb{H}^2$ 
representing the conjugacy class of $\alpha$. Let
$$Y=Y_\alpha=\left(\overline{\mb{H}}^2\setminus\mbox{Fix}(g_\alpha)\right)/
<g_\alpha>,$$
where $\overline{\mb{H}}^2$ is the closed-disk compactification of  
hyperbolic plane $\mb{H}^2$. See Figure \ref{fig;cc2}, where 
$\{a,b\}=\mbox{Fix}(g_\alpha)$. Notice that $Y$ is a 
closed annulus and a neighborhood of $\alpha$ lifts homeomorphically into $Y$.
Suppose that $\beta$ is an essential simple closed curve in $\Sigma$ such that 
$i(\alpha,\beta)\neq 0$. Notice that any lift of $\alpha$ into $Y$ does not 
share endpoints with a lift of $\beta$. Therefore any lift of $\beta$
extends to a properly embedded arc in $Y$. We write $\mb{lift}_\alpha(\beta)$
to denote the set of lifts of $\beta$ into $Y$ which connect the two boundaries
of $Y$ as elements in $\A_0(Y)$. Let $\gamma$ be another essential simple 
closed curve in $\Sigma$ with $i(\alpha,\gamma)\neq 0$. If $\beta$ and $\gamma$
are different, then $b$ and $c$ do not share endpoints for all 
$b\in\mb{lift}_\alpha(\beta)$ and $c\in\mb{lift}_\alpha(\gamma)$. 
The relative twist number is defined by  
$$\tau_\alpha(\beta,\gamma)=\left\{b\cdot c\mid b\in\mb{lift}_\alpha(\beta)\
\mbox{and }c\in\mb{lift}_\alpha(\gamma)\right\}.$$

\indent
From eq. (\ref{eq;dot}), we have $\mbox{diam}(\tau_\alpha(\beta,\gamma))
\leq 2$. In \cite{FLM}, Farb, Lubotzky and Minsky used the following equations 
to prove Theorem \ref{thm;FLM}.
\begin{enumerate}
\item[\textbf{(i)}]
If $t=T_\alpha$ is the leftward Dehn twist on $\alpha$, then 
$\tau_\alpha(\beta,t^n(\beta))\subset\{n,n+1\}$.
\item[\textbf{(ii)}]
If $\beta$ and $\gamma$ intersect $\alpha$, then their geometric intersection 
number bounds their relative twisting, i.e. 
$$\max|\tau_\alpha(\beta,\gamma)|\leq i(\beta,\gamma)+1.$$
\item[\textbf{(iii)}]
If $\beta$, $\gamma$ and $\delta$ intersect $\alpha$, then  
\begin{eqnarray*}
\max\tau_\alpha(\beta,\delta)\leq \max\tau_\alpha(\beta,\gamma)
+\max\tau_\alpha(\gamma,\delta)+2 \\
\min\tau_\alpha(\beta,\delta)\geq \min\tau_\alpha(\beta,\gamma)
+\min\tau_\alpha(\gamma,\delta)-2. 
\end{eqnarray*}
\end{enumerate}

\subsection{The theorem of Masur and Minsky, and of Klarreich}

\begin{figure}[ht]
\begin{center}
\psfrag{0}{\tiny{0}}
\psfrag{1/3}{\tiny{1/3}}
\psfrag{1/2}{\tiny{1/2}}
\psfrag{2/3}{\tiny{2/3}}
\psfrag{1}{\tiny{1}}
\psfrag{3/2}{\tiny{3/2}}
\psfrag{2}{\tiny{2}}
\psfrag{3}{\tiny{3}}
\psfrag{1/0}{\tiny{1/0}}
\psfrag{-1/3}{\tiny{-1/3}}
\psfrag{-1/2}{\tiny{-1/2}}
\psfrag{-2/3}{\tiny{-2/3}}
\psfrag{-1}{\tiny{-1}}
\psfrag{-3/2}{\tiny{-3/2}}
\psfrag{-2}{\tiny{-2}}
\psfrag{-3}{\tiny{-3}}
\includegraphics[width=2.5in,height=2.5in]{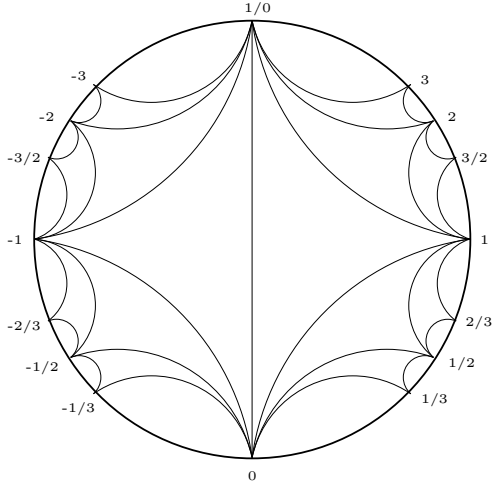}
\end{center}
\caption{Complex of curves of torus}
\label{fig;cc3}
\end{figure}

In this section we study a theorem of Masur and Minsky, and of Klarreich. 
References are \cite{bowditch,hamenstadt,klarreich,MM} and \cite{m-geom}.
If $\Sigma$ is a torus, once-punctured torus or 4-times punctured sphere, then 
the complex of curves is the Farey graph (see Figure \ref{fig;cc3}). 
Minsky showed that $(\C(\Sigma),d_\C)$ is $3\over 2$-hyperbolic space for these
cases (see \cite{m-geom} \S 3 for the proof). For $\xi(\Sigma)=3g+n>4$, the 
following lemma gives an upper bound of $d_\C$. 

\begin{lem}[Lemma 1.1 of \cite{bowditch}]
If $\xi(\Sigma)>4$ then $d_\C(\alpha,\beta)\leq i(\alpha,\beta)+1$.
\end{lem} 
\noindent
Recall the classifications of elements of the mapping class group into 
periodic, reducible and pseudo Anosov elements. An element $h\in\MCG(\Sigma)$ 
is called {\em pseudo Anosov} if there exist $r>1$ and a pair of measured 
foliations $\mc{F}^s$ and $\mc{F}^u$ such that
$$h(\mc{F}^s)={1\over r}\mc{F}^s\quad\mbox{and}\quad h(\mc{F}^u)=r\mc{F}^u.$$ 
If $(\C(\Sigma),d_\C)$ is a bounded metric space with upper bound $K$, then it 
is trivially a $K$-hyperbolic space. 
But from the following proposition, it is clear that 
$\mbox{diam}(\C(\Sigma))=\infty$.
\begin{prop}[Proposition 4.6 of \cite{MM1}]
If $\xi(\Sigma)>4$ then there exists $c>0$ such that, for any pseudo-Anosov
$h\in\MCG(\Sigma)$, $\gamma\in\mc{S}(\Sigma)$ and $n\in\Z$, we have
$$d_\C(h^n(\gamma),\gamma)\geq c|n|.$$
\end{prop}

\indent
In \cite{MM1}, Masur and Minsky proved that $(\C(\Sigma),d_\C)$ is 
$\delta$-hyperbolic for the case $\xi(\Sigma)=3g+n>4$, too. 
\begin{thm}[Theorem 1.1 of \cite{MM}]\label{thm;MM}
$\C(\Sigma_{g,n})$ is a $\delta$-hyperbolic space, where $\delta$ depends only 
on $g$ and $n$.
\end{thm} 
\noindent
Masur and Minsky used Teichm\"uller theory in the proof of Theorem 
\ref{thm;MM}. In their proof, the constant $\delta$ is not constructive because
it contains a compactness arguments on the spaces of quadratic differentials.
In 2002, Bowditch proved Theorem \ref{thm;MM} in a more combinatorial way,
and showed that the number $\delta$ is bounded by a logarithmic function
of $3g+n-4$.    

\indent
Since $\C(\Sigma)$ is $\delta$-hyperbolic, we can consider its Gromov boundary.
For the case of torus, the boundary can be identified with the set of 
irrational numbers. For general cases, it is clear that a sequence of curves 
$\alpha_n$, whose distance from a fixed curve is going to infinity, must 
converge to a maximal lamination. In fact, Klarreich showed that the Gromov 
boundary of complex of curves is homeomorphic to the space of topological
equivalence classes of filling lamination.  
\begin{thm;klarreich}[Klarreich \cite{klarreich}]
There is a homeomorphism $$k:\BC(\Sigma)\to\EL(\Sigma)$$ such that for any 
sequence $\alpha_n$ in $\mc{S}(\Sigma)$, $\alpha_n$ converges to 
$\alpha\in\BC(\Sigma)$ if and only if $\alpha_n$, considered as a subset of 
$\UML(\Sigma)$, converges to $k(\alpha)$.
\end{thm;klarreich}

\indent
Klarreich used Teichm\"uller theory and the results of Masur and Minsky in
\cite{MM}, to prove Theorem \ref{thm;klarreich}. It is clear that 
$\BC(\Sigma)$ is homeomorphic to the Gromov boundary of its 1-skeleton 
$\C_1(\Sigma)$ because they are quasi-isometric. 
Suppose that $\epsilon>0$ satisfies the collar lemma. 
For each $\alpha\in\C_0(\Sigma)$, let
$$T(\alpha)=\{\s\in\T(\Sigma)\mid\ell_\s(\alpha)<\epsilon\}.$$
Then a collection of sets $\T(\alpha_1),\cdots,\T(\alpha_n)$ has nonempty 
intersection if and only if $\alpha_1,\cdots,\alpha_n$ form a simplex in 
$\C(S)$. The set $\T_{el}(\Sigma)$ is defined from  $\T(\Sigma)$, by adding a 
new point $P_\alpha$ for each set $T(\alpha)$ and an interval of length 
$1\over 2$ from $P_\alpha$ to each point in $T(\alpha)$. 
$\T_{el}(\Sigma)$ equipped with the minimal path-metric is called the 
{\em relative Teichm\"uller space} following the terminology of 
Farb(\cite{farb}).
In \cite{MM}, Masur and Minsky showed that $\T_{el}(\Sigma)$ is quasi-isometric
to $\C_1(\Sigma)$, and Klarreich showed that the Gromov boundary of
$\T_{el}(\Sigma)$ is homeomorphic to the space of topological equivalence 
classes of minimal singular foliations on $\Sigma$, which is homeomorphic to 
$\EL(\Sigma)$.

\section{Proof of Main Theorem (Theorem \ref{mainthm})}
\label{chap;pfmainthm}

Suppose that $\s_i\in\T(S)$ converges to $[\lam]\in\PML(S)$ in $\TC(S)$. 
Recall the map 
$$\Phi(\s_i)=\mbox{a pants decomposition whose total length is bounded by }L,$$
where $L$ is a fixed Bers constant and all pants curves are geodesics in 
$\s_i$. Recall also the quotient map $u:\PML(S)\to\UML(S)$ by forgetting 
measure. Notice that $\Phi(\s_i)$ can be considered as a sequence in $\GL(S)$. 
Since $\GL(S)$ is compact, it has a convergent subsequence. 
In this chapter we prove our main theorem.  

\begin{mainthm}
If $\Phi(\s_i)$ converge to $\nu\in\GL(S)$ in Hausdorff metric topology, then 
$u([\lam])\subset\nu$. 
\end{mainthm}

\noindent
Consider the decomposition $u([\lam])=\lam_1\cup\lam_2\cup\cdots\lam_m$ as a 
finite disjoint union of minimal laminations. To prove $\lam_j\subset\nu$ 
for all $1\leq j\leq m$, it is enough to show that $\lam_1\subset\nu$ and 
$\lam_2\subset\nu$, assuming that $\lam_1$ is a simple closed curve and 
$\lam_2$ is not a simple closed curve.

\subsection{Proof of $\lam_1\subset\nu$}\label{sec;lam1}

Recall that $\lam_1$ is a simple closed curve. If $\lam_1\subset\Phi(\s_i)$ for
infinitely many $i$, then it is clear that $\lam_1\subset\nu$. 
If $\lam_1\subset\Phi(\s_i)$ for only finitely many $i$, then there exists 
$N_1>0$ such that $\lam_1\not\subset\Phi(\s_i)$ for all $i>N_1$.
Since $\Phi(\s_i)$ is a pants decomposition, there exists a pants curve 
$\alpha_i$ in $\Phi(\s_i)$ such that $\alpha_i\cap\lam_1\neq\emptyset$ for all
$i>N_1$. Choose $x_i\in\alpha_i\cap\lam_1$ and a limit point $x$ of $x_i$.
By Lemma \ref{lem;gc}, there exists a leaf $\el$ of $\nu$ such that 
$x\in\el$. If $\el=\lam_1$, we are done. 

\begin{figure}[ht]
\begin{center}
\psfrag{r}{$\lam_1$}
\psfrag{l}{$\el$}
\psfrag{x}{$x$}
\includegraphics[width=2in,height=2in]{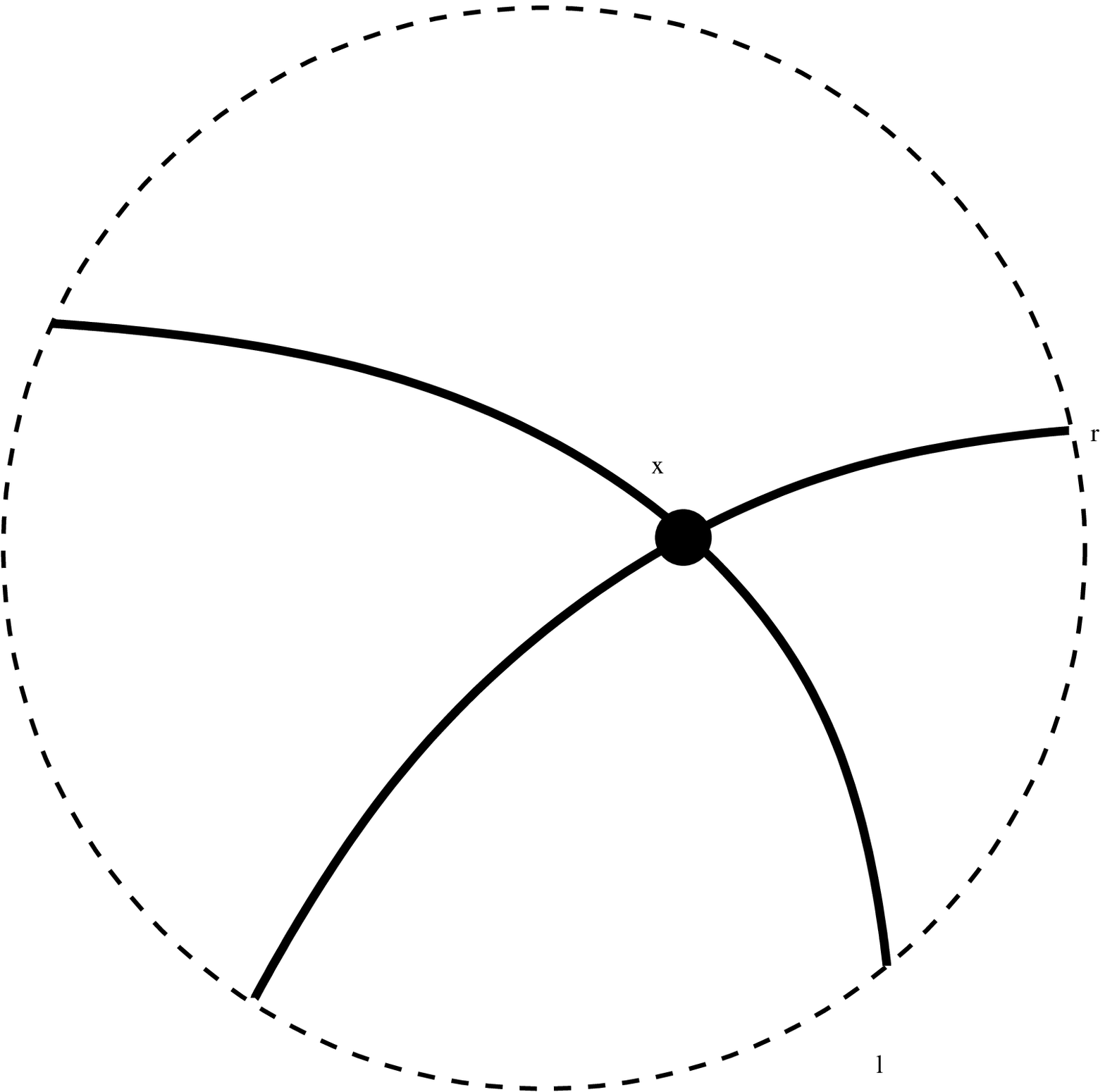}
\end{center}
\caption{The open neighborhood $U$}
\label{fig;pfm1}
\end{figure}

To get a contradiction, suppose that $\el\neq\lam_1$. Choose an open 
neighborhood $U$ of $x$ which is isometric to an open subset of $\mb{H}^2$ 
(see Figure \ref{fig;pfm1}).

\noindent
Since $S$ is compact, by Lemma \ref{lem;HCtopology}, there exists a pants curve
$\alpha_i$ in $\Phi(\s_i)$ such that an arc $\beta_i\subset\alpha_i$ 
approaches to $\el\cap U$. Therefore there exists $N_2>0$ such that 
$$i(\alpha_i,\lam)\geq \int_{\beta_i}d\lam_1=r>0\quad\mbox{for all }i>N_2,$$
where $r>0$ is the transverse measure on $\lam_1$. 

\indent
Choose $\mu_j\in\ML(S)$ which converges to $\lam$ in $\PML(S)$ as in 
Theorem \ref{thm;useful}. Since $\mu_j\to\infty$, there exists a sequence 
$c_j>0$ such that $c_j\mu_j$ converges to $\lam$ in $\ML(S)$ with 
$\lim_{j\to\infty}c_j=0$. Notice that there exists $N_3>0$ which does not 
depend on $i$ such that
$$i(\alpha_i,c_j\mu_j)\geq {r\over 2}\quad \mbox{for all }i,j>N_3.$$
Therefore $i(\alpha_i,\mu_j)$ approaches to $\infty$ as $i,j\to\infty$.
But from eq. (\ref{eq;papa}), we have 
$i(\alpha_i,\mu_i)\leq\el_{\s_i}(\alpha_i)\leq L$ for all $i$.
This is a contradiction.

\subsection{Proof of $\lam_2\subset\nu$}

Recall that $\lam_2$ is not a simple closed curve. Let $F$ be the essential 
subsurface of $S$ which is filled by $\lam_2$. Since $\Phi(\s_i)$ is a pants
decomposition for all $i$, we have $\nu\cap F\neq\emptyset$. 
To get contradictions, suppose that $\lam_2\not\subset\nu$ in the next two
paragraphs. 

\indent
Suppose that $\nu\cap\lam_2=\emptyset$. 
Notice that $\nu\cup\lam_2$ is a geodesic lamination, too. 
By Theorem \ref{thm;SGL}, we can decompose $\nu\cup\lam_2$ as a finite 
disjoint union of minimal lamination, including $\lam_2$, and finite number of 
infinite isolated leaves. Since $\lam_2$ is a filling lamination in $F$, 
any isolated infinite leaf can not intersect $F$. Thus $\nu\cap F=\emptyset$. 
This is a contradiction. 

\begin{figure}[ht]
\begin{center}
\psfrag{r}{$\lam_2$}
\psfrag{l}{$\el$}
\includegraphics[width=2in,height=2in]{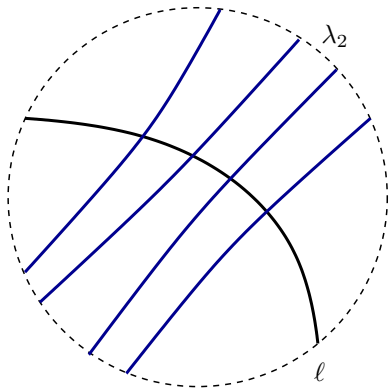}
\end{center}
\caption{The open neighborhood $V$}
\label{fig;pfm2}
\end{figure}

\indent
Suppose that $\nu\cap\lam_2\neq\emptyset$. There exists a leaf $\ell$ of
$\nu$ which intersects $\lam_2$ transversely. Choose an open neighborhood $V$ 
which is isometric to an open subset of $\mb{H}^2$, and in which $\el$ 
intersects $\lam_2$ transversely (see Figure \ref{fig;pfm2}).
As in Section \ref{sec;lam1}, there exists a pants curve 
$\alpha_i$ in $\Phi(\s_i)$ and arc $\beta_i\subset\alpha_i$ such that
$\beta_i$ converges to $\el\cap V$. Therefore there exist $r>0$ and $N>0$ 
such that 
$$i(\alpha_i,\lam)\geq\int_{\beta_i}d\lam_2=r>0\quad\mbox{for all }i>N.$$
As in Section \ref{sec;lam1}, using Theorem \ref{thm;useful}, we can show 
that this is a contradiction.

\section{Proof of Corollary \ref{cor;twist} and Theorem \ref{thm;surgery}}
\label{chap;ct}

In this chapter we prove Corollary \ref{cor;twist} and Theorem 
\ref{thm;surgery}, and show that Theorem \ref{thm;filling} comes from
Theorem \ref{thm;surgery}.
 
\subsection {Proof of Corollary \ref{cor;twist}}\label{sec;twist} 

Suppose that $\s_i\in\T(S)$ converges to $[\lam]\in\PML(S)$ in $\TC(S)$, and 
let $u([\lam])= \lam_1\cup\lam_2\cup\cdots\cup\lam_m$ be the decomposition as a
finite disjoint union of minimal laminations. 

\begin{figure}[ht]
\begin{center}
\psfrag{U}{$U$}
\psfrag{r}{$\lam_1$}
\psfrag{Y}{$Y$}
\psfrag{S}{$S$}
\psfrag{p}{$p$}
\includegraphics[width=2.5in,height=2in]{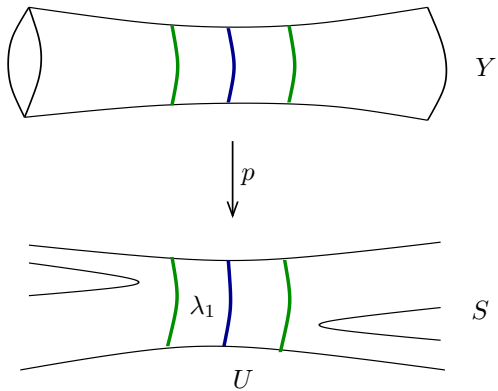}
\end{center}
\caption{The annular cover $Y$}
\label{fig;ct1}
\end{figure}

\noindent
Suppose that $\lam_1$ is a simple
closed curve. Suppose also that $\alpha_i$ is a pants curve in $\Phi(\s_i)$ and
$\alpha_i\cap\lam_1\neq\emptyset$ for all $i$.  
Construct an annular covering $Y$ of $S$ in which a neighborhood $U$ of 
$\lam_1$ lifts homeomorphically (see Figure \ref{fig;ct1}). 
We may assume that $U$ is a closed collar around $\lam_1$ and $U$ 
does not intersect $\lam_j$ for all $j\neq 1$. 
Let $a_i\in\mb{lift}(\alpha_i)$. In this section we prove Corollary 
\ref{cor;twist}.

\begin{twistthm}
$|a_1\cdot a_i|$ approaches to infinity as $i$ increases.
\end{twistthm} 

\begin{proof}
To get a contradiction, suppose that $|a_1\cdot a_i|$ does not approach to 
$\infty$. Then there exists a subsequence of $a_i$, which we will call $a_i$ 
again for the sake of simplicity, such that $|a_1\cdot a_i|=k$ for all $i$ 
for some $k\in\N\cup\{0\}$. We may assume that $a_1\cdot a_i=k$ without loss
of generality. 

\begin{figure}[ht]
\begin{center}
\psfrag{H}{${\overline{\mb{H}}}^2$}
\psfrag{p}{$p$}
\psfrag{U}{$U$}
\psfrag{K}{$K$}
\psfrag{b}{$\beta_i$}
\includegraphics[width=1.8in,height=3in]{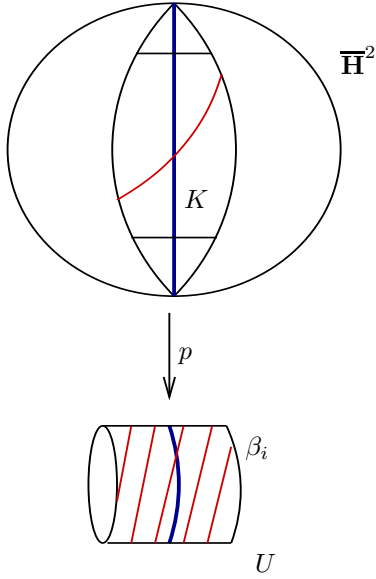}
\end{center}
\caption{A lift of $\beta_i$ in $K$}
\label{fig;ct2}
\end{figure}

\noindent
Suppose that $i\geq 3$. From eq. (\ref{eq;dot}), we have
$a_1\cdot a_i=a_1\cdot a_2+a_2\cdot a_i+\Delta$, where $\Delta\in\{-1,0,1\}$. 
Therefore 
$$ |a_2\cdot a_i|\leq 1\quad\mbox{for all }i\geq 3.$$
This equation implies that $\alpha_i=p(a_i)$ does not change so much in $U$  
for $i\geq 2$. Let $\beta_i$ be a component of $\alpha_i\cap U$ (see Figure 
\ref{fig;ct2}). 
We can find a compact set $K\subset {\mb{H}}^2$ such that there 
exists a lift of $\beta_i$ in $K$ for all $i$. A geodesic arc in $\mb{H}^2$ is 
determined by its two endpoints. Therefore there exists a subsequence 
of $\beta_i$, which we will call $\beta_i$ again, which converges to 
a geodesic arc $\beta$ in Hausdorff metric topology with 
$i(\beta,\lam_1)\neq 0$. For this subsequence $$\beta_i\subset\alpha_i
\subset\Phi(\s_i),$$  
$\Phi(\s_i)$ still converges to $\nu$ in Hausdorff metric topology.
Notice that $\lam_1\not\subset\nu$. This is a contradiction to Theorem 
\ref{mainthm}.
\end{proof}

\subsection{Proof of Theorem \ref{thm;surgery}}\label{sec;surgery}

\begin{figure}[ht]
\begin{center}
\psfrag{b}{$\beta_i$}
\psfrag{F}{$F$}
\includegraphics[width=3in,height=2in]{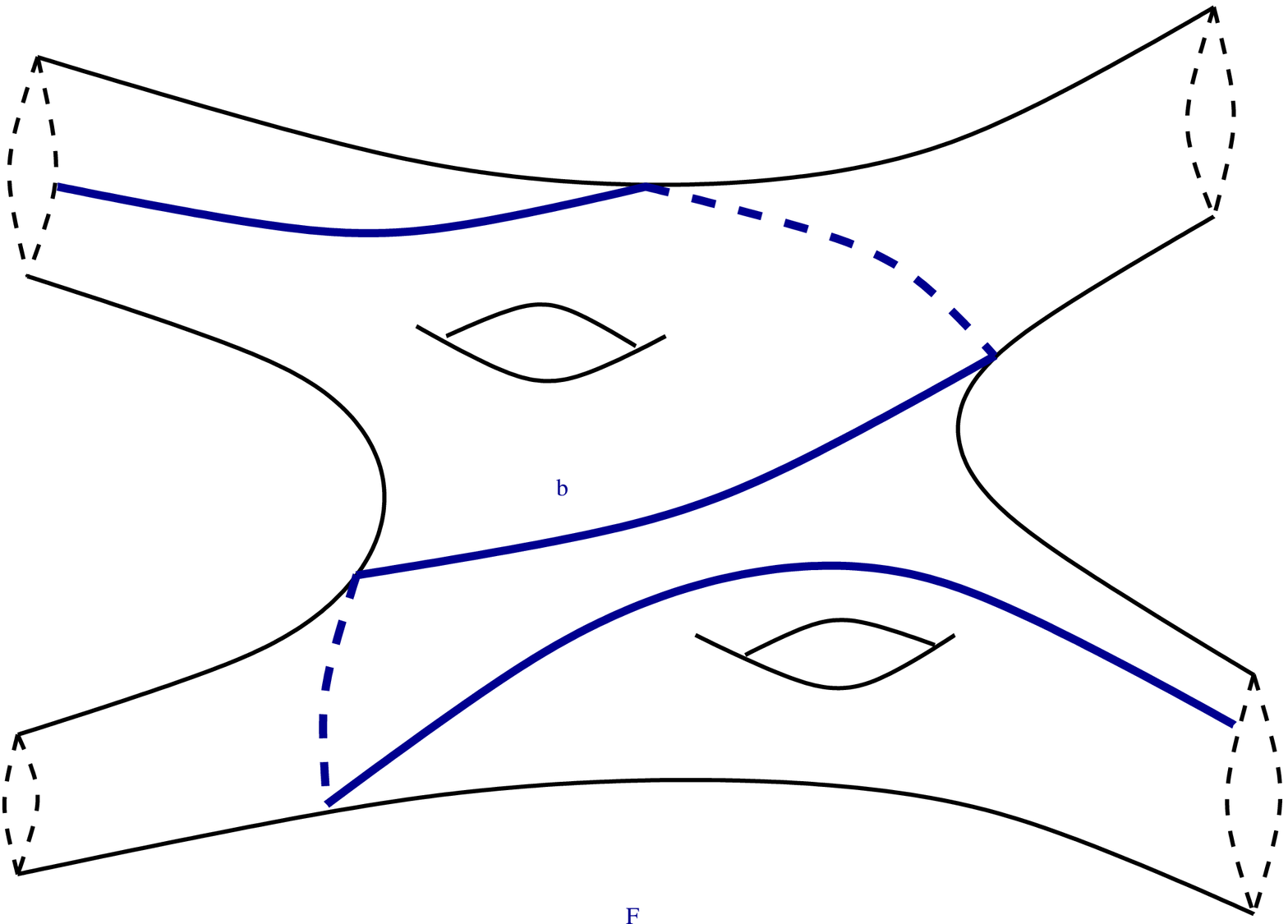}
\end{center}
\caption{$\beta_i$}
\label{fig;ct3}
\end{figure}

In this section we prove Theorem \ref{thm;surgery}.
Suppose that $\lam_2$ is not a simple closed curve and let $F$ be the 
subsurface of $S$ which is filled by $\lam_2$ (see Figure \ref{fig;ct3}). 

\noindent
Notice that there exists a pants curve $\alpha_i$ in $\Phi(\s_i)$ with 
$\alpha_i\cap F\neq\emptyset$ for all $i$. 
Let $\beta_i$ be a component of $\alpha_i\cap F$. 
As in Lemma 2.2 of \cite{MM2}, let
\begin{eqnarray*}
&&\Tilde{\beta_i}\mbox{ be a non-peripheral essential component of boundary}\\
&&\quad\quad\mbox{of regular neighborhood of }\beta_i\cup\partial\overline{F},
\end{eqnarray*}

\noindent
where $\overline{F}$ is the completion of $F$ with path-metric. Note that two 
different components of $\partial\overline{F}$ could be a same curve in $S$.

\begin{figure}[ht]
\begin{center}
\psfrag{b}{$\beta_i$}
\psfrag{1}{\textbf{Case I}}
\psfrag{2}{\textbf{Case II}}
\psfrag{B}{$\partial\overline{F}$}
\includegraphics[width=4.5in,height=2in]{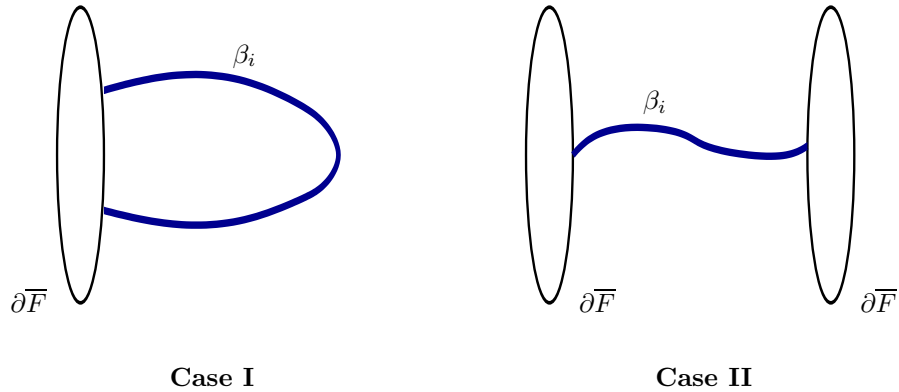}
\end{center}
\caption{Two Cases}
\label{fig;ct4}
\end{figure}

\indent
Notice that if $\beta_i$ is a closed curve, then $\Tilde{\beta_i}$ is homotopic
to $\beta_i$. Notice also that if $\beta_i$ is an arc, then there are two cases
as in Figure \ref{fig;ct4}.

\noindent
Recall that $\lam_2$ is not a closed curve. Since $F$ is filled by $\lam_2$, 
it can not be a disk, an annulus or a pants. Therefore for both cases, the 
regular neighborhood of $\beta_i\cup\partial\overline{F}$ has a boundary 
component which is non-peripheral and essential in $F$.   
Suppose that $\partial\overline{F}=\{\gamma_1,\cdots,\gamma_k\}$ and let
$$\el_{\s_i}(\partial F)=\el_{\s_i}(\gamma_1)+\cdots\el_{\s_i}(\gamma_k).$$
Let $\el_{\s_i}(\Tilde{\beta}_i)$ be the length of geodesic representative of 
$\Tilde{\beta}_i$ in $\s_i$. From the definition of $\Tilde{\beta}_i$, we have 
\begin{equation}\label{eq;simple}
\el_{\s_i}(\Tilde{\beta}_i)<2\el_{\s_i}(\alpha_i)+\el_{\s_i}(\partial F).
\end{equation} 

\indent
We now prove Theorem \ref{thm;surgery}. 

\begin{thm;surgery}
The geodesic representative of $\Tilde{\beta_i}$ converges to $\lam_2$ in 
$\UML(F)$.   
\end{thm;surgery} 

\noindent
To prove Theorem \ref{thm;surgery}, it is enough to show that the geodesic
representative of $\Tilde{\beta_i}$ converges to $\lam_2$ in $\UML(S)$. 
For any limit point $[\beta]$ of $[\Tilde{\beta_i}]$ in $\PML(S)$, we will show
that $i(\beta,\lam_2)=0$. Then Theorem \ref{thm;surgery} follows from the fact
that $\lam_2$ fills $F$.

\begin{lem}
If $[\beta]$ is a limit point of $[\Tilde{\beta_i}]$ in $\PML(S)$, then 
$i(\beta,\lam_2)=0$.
\end{lem}

\begin{proof}
Suppose that $[\beta]$ is a limit point of $[\Tilde{\beta_i}]$. After possibly
restricting to a subsequence, we may assume that $[\Tilde{\beta_i}]$ converges 
to $[\beta]$ in $\PML(S)$. There exist a constant $K>0$ and a sequence $b_i>0$ 
such that $b_i\Tilde{\beta_i}$ converges to $\beta$ in $\ML(S)$ with 
$b_i\leq K$ for all $i$. 

\indent
Choose $\mu_i\in\ML(S)$ as in Theorem \ref{thm;useful}. Since $\mu_i$ converges
to $[\lam]\in\PML(S)$ with $\mu_i\to\infty$, there exists a sequence $c_i>0$ 
such that $c_i\mu_i$ converges to $\lam$ in $\ML(S)$ with $c_i\to 0$. 
By eq. (\ref{eq;papa}), there exists $\Gamma>0$ which does not depend on $i$ 
such that 
$i(\partial\overline{F},\mu_i)\leq \el_{\s_i}(\partial\overline{F})\leq 
i(\partial\overline{F},\mu_i)+\Gamma$.
Therefore
$$i(\partial\overline{F},c_i\mu_i)\leq c_i\el_{\s_i}(\partial\overline{F})\leq 
i(\partial\overline{F},c_i\mu_i)+c_i\Gamma.$$
Since $i(\partial\overline{F},\lam)=0$, from the continuity of the intersection
number, we have
\begin{equation}\label{eq;crucial}
\lim_{i\to\infty}c_i\el_{\s_i}(\partial\overline{F})=0.
\end{equation}
From eq. (\ref{eq;papa}) and eq. (\ref{eq;simple}), we have
$$i(\Tilde{\beta_i},\mu_i)\leq \el_{\s_i}(\Tilde{\beta}_i)\leq 
2\el_{\s_i}(\alpha_i)+\el_{\s_i}(\partial\overline{F})$$
Therefore 
$$i(b_i\Tilde{\beta_i},c_i\mu_i)\leq K\left(2c_i L+c_i\el_{\s_i}(
\partial\overline{F})\right).$$ Hence from eq. (\ref{eq;crucial}) and the 
continuity of the intersection number, we have $i(\beta,\lam)=0$. 
Thus $i(\beta,\lam_2)=0$.
\end{proof}

\subsection{Theorem \ref{thm;filling} follows from Theorem \ref{thm;surgery}}

Suppose that $\s_i\in\T(S)$ converges to $[\lam]\in\PML(S)$ in $\TC(S)$, and 
$\alpha_i$ is a pants curve in $\Phi(\s_i)$. If $\lam$ is a filling lamination,
from Theorem \ref{thm;surgery}, then $\alpha_i$ converges to $u([\lam])$ in 
$\UML(S)$. Therefore if we identify $\BC(S)$ with $\EL(S)$ via the 
homeomorphism $k$ in Theorem \ref{thm;klarreich}, we have
\begin{thm;filling} 
If $\lam$ is a filling lamination, then $\alpha_i$ converges to $u([\lam])$ in 
$\C(S)\cup\BC(S)$.
\end{thm;filling}



\begin{thebibliography}{00}

\bibitem {abikoff} W. Abikoff, \textit{Kleinian groups-geometrically finite 
and geometrically perverse}, Geometry of group representation (Boulder, CO, 
1987), 1-50, Contemp. Math. 74, Amer. Math. Soc., 1988.

\bibitem {agol} I. Agol, \textit{Tameness of hyperbolic 3-manifolds}, preprint,
arXiv:math. GT/0405568.

\bibitem {ahlfors} L.V. Ahlfors, \textit{On quasiconformal mappings}, J. 
Analyse Math. {\bf 3}(1954), 1-58.

\bibitem {BP} R. Benedetti and C. Petronio, \textit{Lectures on hyperbolic 
geometry}, Springer-Verlag, 1992.

\bibitem{beardon} A.F. Beardon, \textit{The geometry of discrete group}, 
Springer-Verlag, New York, 1983.

\bibitem {bers} L. Bers, \textit{Spaces of Riemann surfaces as bounded 
domains}, Bull. Amer. Math. Soc. {\bf 66}(1960), 98-103.

\bibitem {bers1} L. Bers, \textit{On boundary of Teichm\"uller spaces
and on Kleinian groups I}, Ann. of Math. {\bf 91}(1970), 570-600.

\bibitem {bers2} L. Bers, \textit{Spaces of degenerating Riemann surfaces}, 
Discontinuous groups and Riemann surfaces, Ann. of Math. Stud 79, Princeton 
Univ. Press, 1974, pp. 43-59.

\bibitem {bers3} L. Bers, \textit{An inequality for Riemann surfaces}, 
Differential geometry and complex analysis, 87-93, Springer-Verlag, 1985.

\bibitem {bonahon} F. Bonahon, \textit{Bouts des vari\'et\'es hyperboliques 
de dimension 3}, Ann. of Math. {\bf 124}(1986), 71-158.
 
\bibitem {bon} F. Bonahon, \textit{The geometry of Teichm\"uller space via 
geodesic currents}, Invent. Math. {\bf 92}(1988), 139-162.
 
\bibitem {bon2} F. Bonahon, \textit{Transverse H\"oder distributions for 
geodesic laminations}, Topology {\bf 36}(1997), no.1, 103-122.
 
\bibitem {bon3} F. Bonahon, \textit{Geodesic laminations with transverse 
H\"oder distributions}, Ann. Sci. Ecole Norm. Sup.(4) {\bf 30}(1997), 
205-240.

\bibitem {bon4} F. Bonahon, \textit{Geodesic laminations on surfaces}, 
Laminations and foliations in dynamics, geometry and topology (Stony Brook, NY,
1998), 1-37, Contemp. Math. 269, Amer. Math. Soc., 2001.
 
\bibitem {BoS} M. Bonk and O. Schramm, \textit{Embeddings of Gromov hyperbolic 
spaces}, Geom. Funct. Anal. {\bf 10}(2000), no.2, 266-306.

\bibitem {bowditch} B.H. Bowditch, \textit{Intersection numbers and the 
hyperbolicity of the curve complex}, preprint, 2002.

\bibitem {bowers} P.L. Bowers, \textit{Negatively curved graph and planar 
metrics with applications to type}, Michigan Math. J. {\bf 45}(1989), no.1,
31-53.

\bibitem {bridson} M. Bridson, \textit{Geodesics and curvature in metric 
simplicial complexes}, Group theory from a geometrical viewpoint (Triests, 
1990), 373-463, World Sci. Publishing, 1991.

\bibitem {BH} M. Bridson and A. Haefliger, \textit{Metric spaces of 
non-positive curvature}, Springer-Verlag, 1999.

\bibitem {brock} J.F. Brock, \textit{Continuity of Thurston's length function},
Geom. Funct. Anal. {\bf 10}(2000), no.4, 741-797. 

\bibitem {brock2} J.F. Brock, \textit{Boundaries of Teichm\"uller space and 
end-invariants for hyperbolic 3-manifolds}, Duke Math. J. {\bf 106}(2001), 
no.3, 527-552. 

\bibitem {buser} P. Buser, \textit{Geometry and spectra of compact Riemann 
surfaces}, Birkh\"{a} user, 1992.

\bibitem {BS} P. Buser and M. Sepp\"al\"a, \textit{Symmetric pants 
decompositions of Riemann surfaces}, Duke. Math. J. {\bf 67}(1992), no.1, 
39-55.

\bibitem {CG} D. Calegari and D. Gabai, \textit{Shrinkwrapping and the taming 
of hyperbolic 3-manifolds}, preprint, arXiv:math.GT/0407161. 

\bibitem {CEG} R.D. Canary, D.B.A. Epstein and P. Green, \textit{Notes on notes
 of Thurston}, Analytical and geometrical aspects of Hyperbolic spaces (
D.B.A. Epstein, ed.), London Math. Lecture Note Ser. 111, Cambridge Univ. 
Press, 1987, pp. 3-92.

\bibitem {canary} R.D. Canary, \textit{Ends of hyperbolic 3-manifolds}, 
J. Amer. Math. Soc. {\bf 6}(1993), no.1, 1-35.

\bibitem {cannon} J. Cannon, \textit{The theory of negatively curved spaces and
groups}, Ergodic theory, symbolic dynamics, and hyperbolic spaces 
(Trieste, 1989), Oxford Univ. Press, 1991, 315-369.

\bibitem {CB} A.J. Casson and S.A. Bleiler, \textit{Automorphisms of surfaces 
after Nielsen and Thurston}, Cambridge University Press, 1988.

\bibitem {farb} B. Farb, \textit{Relative hyperbolic and automatic groups with 
applications to negatively curved manifolds}, Ph.D. thesis, Princeton Univ., 
1994. 

\bibitem {FLM} B. Farb, A. Lubotzky and Y. Minsky, \textit{Rank-1 phenomena for
mapping class group},  Duke. Math. J. {\bf 106}(2001), no.3, 581-597.

\bibitem {FLP} A. Fathi, F. Laudenbach and V. Poenaru, \textit{Travaux de 
Thurston sur les surfaces}, vol. 66-67, Asterisque, 1979.

\bibitem {gabai} D. Gabai, \textit{3 lectures on foliations and laminations on 
3-manifolds}, Laminations and foliations in dynamics, geometry and topology
(Stony Brook, NY, 1998), 87-109, Contemp. Math. 269, Amer. Math. Soc., 2001.

\bibitem {gardiner} F.P. Gardiner, \textit{A correspondence between laminations
 and quadratic differentials}, Complex Variables {\bf 6}(1986), 363-375. 

\bibitem {GM} F.P. Gardiner and H. Masur, \textit{Extremal length geometry 
of Teichm\"uller space}, Complex Variables, {\bf 16}(2000), 209-237. 

\bibitem {gersten} S.M. Gersten, \textit{Introduction to hyperbolic and 
automatic group}, Summer school in group theory in Banff, 1996, 45-70, CRM 
Proc. Lecture Notes 17, Amer. Math. Soc., 1999. 

\bibitem {gromov} H. Gromov, \textit{Hyperbolic groups}, Essays in Group Theory
(S.M. Gersten, editor), MSRI Publications no.8, Springer-Verlag, 1978.

\bibitem {hamenstadt} U. Hamenst\"{a}dt, \textit{Train tracks and the Gromov 
boundary of the complex of curves}, preprint, arXiv:math.GT/0409611. 

\bibitem {harer} J.L. Harer, \textit{The virtual cohomological dimension of
the mapping class group of an orientable surface}, Invent. Math. 
{\bf 84}(1986), 157-176. 

\bibitem {harvey} W.J. Harvey, \textit{Boundary structure of the modular 
group}, Riemann surfaces and related topics, 245-251, Ann. of Math. Stud. 97, 
Princeton Univ. Press, 1981.

\bibitem {hatcher} A.E. Hatcher, \textit{Measured lamination spaces for 
surfaces, from the topological viewpoint}, Topology Appl. {\bf 30}(1988),
63-88. 

\bibitem {HT} A. Hatcher and W. Thurston, \textit{A presentation for the 
mapping class group of a closed orientable surface}, Topology {\bf 19}(1980),
221-237. 

\bibitem {HM} J. Hubbard and H. Masur, \textit{Quadratic differentials and 
foliations}, Acta Math. {\bf 142}(1979), no.3-4, 221-274.

\bibitem {IT} Y. Imayoshi and M. Taniguchi, \textit{An introduction to 
Teichm\"uller spaces}, Springer-Verlag, 1992.

\bibitem {ivanov} N.V. Ivanov, \textit{Automorphisms of complexes of curves 
and of Teichm\"uller spaces}, Internat. Math. Res. Notices {\bf 1997}, no.14,
651-666.

\bibitem {KM} V.A. Kaimanovich and H. Masur, \textit{The Poisson boundary of 
the mapping class group}, Invent. Math. {\bf 125}(1996), no.2, 221-264.

\bibitem {keen} L. Keen, \textit{Collars on Riemann surfaces}, Discontinuous
groups and Riemann surfaces, Ann. of Math. Studies, No.79, Princeton Univ. 
Press, 1974. 

\bibitem {kerckhoff} S.P. Kerckhoff, \textit{The asymptotic geometry of 
Teichm\"uller space}, Topology {\bf 19}(1980), 23-41.

\bibitem {kerckhoff2} S.P. Kerckhoff, \textit{The Nielsen realization problem},
Ann. Math. {\bf 117}(1983), 235-265.

\bibitem {KT} S.P. Kerckhoff and W.P. Thurston, \textit{Non-continuity of the 
action of the modular group at Bers' boundary of Teichm\"uller space},
Invent. Math. {\bf 100}(1990), 25-47.

\bibitem {kim} Y.D. Kim, \textit{A theorem on discrete, torsion free subgroups 
of $Isom\, H^n$}, Geom. Dedicata {\bf 109}(2004), 51-57.

\bibitem {klarreich} E. Klarreich, \textit{The boundary at infinity of the 
curve complex and the relative Teichm\"uller space}, preprint, 1998.

\bibitem {korkmaz} M. Korkmaz, \textit{Automorphisms of complexex of curves on 
punctured spheres and on punctured tori}, Topology and its Applications 
{\bf 95}(1999), 85-111.

\bibitem {kra} I. Kra, \textit{Quadratic differentials}, Rev. Roumaine Math. 
Pures Appl. {\bf 39}(1994), no.8, 751-787.

\bibitem {krushkal} S.L. Krushkal, \textit{Teichm\"uller spaces are not 
starlike}, Ann. Acad. Sci. Fenn. Ser.A/Math. {\bf 20}(1995), no.1, 167-173.

\bibitem {KS} R.S. Kulkarni and P.B. Shallen, \textit{On Ahlfors' finiteness
theorem}, Adv. Math. {\bf 76}(1989), no.2, 155-169.

\bibitem {levitt} G. Levitt, \textit{Foliations and laminations on hyperbolic 
surfaces}, Topology {\bf 22}(1983), no.2, 119-135. 

\bibitem {luo} F. Luo, \textit{Grothendieck's reconstruction principle and
2-dimensional topology and geometry}, Commun. Contemp. Math. {\bf 1}(1999), 
no.2, 125-153.

\bibitem {luo2} F. Luo, \textit{Geodesic length functions and Teichm\"uller 
space}, Electron. Res. Announc. Amer. Math. Soc. {\bf 2}(1996), no.1, 34-41.

\bibitem {luo3} F. Luo, \textit{Automorphisms of the complex of curves}, 
Topology {\bf 39}(2000), 283-298.

\bibitem {marden} A. Marden, \textit{The geometry of finitely generated
Kleinian groups}, Ann. of Math. {\bf 99}(1974), 383-462. 

\bibitem {maskit} B. Maskit, \textit{Comparison of hyperbolic and extremal 
lengths}, Ann. Acad. Sci. Fenn. Series A. I. Mathematica
{\bf 10}(1985), 381-386.

\bibitem {maskit2} B. Maskit, \textit{On boundary of Teichm\"uller spaces
and on Kleinian groups II}, Ann. of Math. {\bf 91}(1970), 607-639.

\bibitem {maskit3} B. Maskit, \textit{A picture of moduli space}, Invent. Math.
 {\bf 126}(1996), 341-390.

\bibitem {maskit4} B. Maskit, \textit{Kleinian groups}, Springer-Verlag, 1988.

\bibitem {masur} H. Masur, \textit{On a class of geodesics in Teichm\"uller 
space}, Ann. of Math. {\bf 102}(1975), 205--221. 

\bibitem {masur1} H. Masur, \textit{Two boundaries of Teichm\"uller space}, 
Duke Math. J. {\bf 49}(1982), no.1, 183-190. 

\bibitem {MM} H. Masur and Y. Minsky, \textit{Geometry of the complex of
curves I, hyperbolicity}, Invent. Math. {\bf 138}(1999), 103-149.

\bibitem {MM1} H. Masur and Y. Minsky, \textit{Unstable quasi-geodesics
in Teichm\"uller space}, Contemp. Math. 256, Amer. Math. Soc., 2000.

\bibitem {MM2} H. Masur and Y. Minsky, \textit{Geometry of the complex of 
curves II, hierarchical structure}, Geom. Funct. Anal. {\bf 10}(2000), 902-974.

\bibitem {MW} H. Masur and M. Wolf, \textit{Teichm\"uller space is not Gromov 
hyperbolic}, Ann. Acad. Sci. Fenn. Ser.A/Math. {\bf 20}(1995), no.2, 259-267.

\bibitem {MP} J. McCarthy and A. Papadopoulos, \textit{Dynamics on Thurston's 
sphere of projective measured foliations}, Comment. Math. Helvetici 
{\bf 64}(1989), 133-166.

\bibitem {MP1} J. McCarthy and A. Papadopoulos, \textit{The visual sphere of 
 Teichm\"uller space and a theorem of Masur-Wolf}, Ann. Acad. Sci. Fenn. Math.
 {\bf 24}(1999), 147-154.

\bibitem {mccullough} D. McCullough, \textit{Compact submanifolds of 
3-manifolds with boundary}, Quart. J. Math. Oxford {\bf 37}(1986), 299-306.

\bibitem {mcmullen} C. McMullen, \textit{Cusps are dense}, Ann. of Math. 
{\bf 133}(1991), 217-247.

\bibitem {mosher} L. Mosher, \textit{Stable Teichm\"uller quasigeodesics 
and ending laminations}, Geom. Topol. {\bf 7}(2003), 33-90.

\bibitem {m-thesis} Y. Minsky, \textit{Harmonic maps, length, and energy in
Teichm\"uller space}, J. Differential Geom. {\bf 35}(1992), 151-217.

\bibitem {m-teich} Y. Minsky, \textit{Teichm\"uller geodesics and ends
of hyperbolic 3-manifolds}, Topology {\bf 32}(1993), 625-647.

\bibitem {m-geom} Y. Minsky, \textit{A geometric approach to the complex of 
curves on a surface}, Topology and Teichm\"uller space (Katinkulta, 1995), 
149-158, World Sci. Publishing, 1996.

\bibitem {m-extremal} Y. Minsky, \textit{Extremal length estimate and product 
regions in Teichm\"uller space}, Duke. Math. J. {\bf 83}(1996), 249-286.

\bibitem {m-quasi} Y. Minsky, \textit{Quasi-projections in Teichm\"uller 
space}, J. Reine Angew. Math. {\bf 473}(1996), 121-136.

\bibitem {m-torus} Y. Minsky, \textit{The classification of punctured-torus
group}, Ann. of Math. {\bf 149}(1999), 559-626.

\bibitem {m-klein} Y. Minsky, \textit{Kleinian groups and the complex of 
curves}, Geometry and Topology {\bf 4}(2000), 117-148.

\bibitem {m-combi} Y. Minsky, \textit{Combinatorial and geometrical aspects of
hyperbolic 3-manifolds}, London Math. Soc. Lecture Note Ser. 299, Cambridge 
Univ. Press, 2003.

\bibitem {m-bounded} Y. Minsky, \textit{Bounded geometry for Kleinian groups}, 
Invent. Math. {\bf 146}(2001), no.1, 143-192.

\bibitem {m-cdm} Y. Minsky, \textit{End invariants and the classification of 
hyperbolic 3-manifolds}, Current developments in mathematics, 2002, 181-217,
Int. Press, 2003. 

\bibitem {m-ending} Y. Minsky, \textit{The classification of Kleinian surface
groups I: Models and bounds}, preprint, arXiv:math.GT/0302208.

\bibitem {m-qconvex} Y. Minsky, \textit{Quasiconvexity in the curve complex}, 
preprint, arXiv:math. GT/0307083.

\bibitem {m-ending1} Y. Minsky, \textit{The classification of Kleinian surface
groups II: The ending lamination conjecture}, preprint, arXiv:math. GT/0412006.

\bibitem {MO} J.W. Morgan and J. Otal, \textit{Relative growth rates of closed 
geodesics on a surface under varying hyperbolic structures}, 
Comment. Math. Helv. {\bf 68}(1993), no.2, 171-208. 

\bibitem {papa} A. Papadopoulos, \textit{On Thurston's boundary of
Teichm\"uller space and the extention of earthquakes}, Topology Appl. 
{\bf 41}(1991), no.3, 147-177.

\bibitem {penner} R.C. Penner, \textit{An introduction to train tracks}, 
Low-dimensional topology and Kleinian groups, 77-90, London Math. Soc. Lecture 
Note Ser. 112, Cambridge Univ. Press, 1986. 

\bibitem {PH} R.C. Penner and J.L. Harer, \textit{Combinatorics of train 
tracks}, Annals of Mathematics Studies 125, Princeton University Press, 
Princeton, 1992.

\bibitem {PM} A. Portolano and I. Maniscalco, \textit{Curves as measured 
foliation on noncompact surface}, Rend. Circ. Mat. Palermo(2) {\bf 42}(1993), 
no.2, 161-180. 

\bibitem {ratcliffe} J.G. Ratcliffe, \textit{ Foundations of hyperbolic 
manifolds}, Springer-Verlag, New York, 1994. 

\bibitem {schmutz} P. Schmutz, \textit{Riemann surfaces with shortest geodesic 
of maximal length}, Geom. Funct. Anal. {\bf 3}(1993), no.6, 564-631.

\bibitem {schmutz2} P. Schmutz, \textit{Systoles on Riemann surfaces}, 
Manuscripta Math. {\bf 85}(1994), no.3-4, 429-447. 

\bibitem {scott} G.P. Scott, \textit{Compact submanifolds of 3-manifolds}, 
J. London Math. Soc. {\bf 7}(1973), 246-250. 

\bibitem {short} Edited by H. Short, \textit{Notes on word hyperbolic groups}, 
Group theory from geometric viewpoint (Trieste, 1990), 3-63, World Sci. 
Publishing, 1991.

\bibitem {sullivan} D. Sullivan, \textit{On the ergodic theory at infinity of
an arbitrary discrete group of hyperbolic motions}, Riemann surfaces 
and related topics: Proceedings of the 1978 Stony Brook Conference, 
Ann. of Math. Stud. {\bf 97}, Princeton, 1981.

\bibitem {thurston} W. Thurston, \textit{3-dimensional geometry and topology}, 
Princeton University Press, 1997, (S. Levy, ed.).

\bibitem {thurston1} W. Thurston, \textit{The geometry and topology of 
3-manifolds}, Princeton University Lecture Notes, 1982.
  
\bibitem {thurston2} W. Thurston, \textit{Three dimensional manifolds, Kleinian
groups and hyperbolic geometry}, Bull. Amer. Math. Soc. {\bf 6}(1982), 357-381.

\bibitem {thurston3} W. Thurston, \textit{Minimal stretch maps between 
hyperbolic surfaces}, preprint, arXiv:math.GT/9801039. 

\bibitem {thurston4} W. Thurston, \textit{Hyperbolic structures on 3-manifolds 
II: Surface groups and 3-manifolds which fiber over the circle}, preprint, 
arXiv:math. GT/9801045. 

\bibitem {thurston5} W. Thurston, \textit{On the geometry and dynamics of 
diffeomorphisms of surfaces}, Bull. Amer. Math. Soc. {\bf 19}(1988), no.2, 
417-431.

\bibitem {weiss} H. Weiss, \textit{The geometry of measured geodesic 
laminations and measured train tracks}, Ergodic Theory Dynam. Systems 
{\bf 9}(1989), no.3, 587-604.

\end{thebibliography}
\end{document}